

\input amstex

\define\scrO{\Cal O}
\define\Pee{{\Bbb P}}
\define\Zee{{\Bbb Z}}
\define\Cee{{\Bbb C}}
\define\Ar{{\Bbb R}}

\define\dbar{\bar {\partial}}
\define\dirac{\rlap{/}\partial}
\define\Pic{\operatorname{Pic}}
\define\Ker{\operatorname{Ker}}
\define\Coker{\operatorname{Coker}}
\define\Sym{\operatorname{Sym}}
\define\Hom{\operatorname{Hom}}
\define\Ext{\operatorname{Ext}}
\define\Id{\operatorname{Id}}
\define\ch{\operatorname{ch}}
\define\Todd{\operatorname{Todd}}

\define\proof{\demo{Proof}}
\define\endproof{\qed\enddemo}
\define\endstatement{\endproclaim}
\define\theorem#1{\proclaim{Theorem #1}}
\define\lemma#1{\proclaim{Lemma #1}}
\define\proposition#1{\proclaim{Proposition #1}}
\define\corollary#1{\proclaim{Corollary #1}}
\define\claim#1{\proclaim{Claim #1}}

\define\section#1{\specialhead #1 \endspecialhead}
\define\ssection#1{\medskip\noindent{\bf #1}}

\documentstyle{amsppt}

\topmatter
\title
Obstruction bundles, semiregularity, and Seiberg-Witten invariants
\endtitle
\rightheadtext{Obstruction bundles, semiregularity, Seiberg-Witten invariants}
\author {Robert Friedman  and John W. Morgan}
\endauthor
\address Department of Mathematics, Columbia University, New York,
NY 10027, USA\endaddress
\email rf\@math.columbia.edu, jm\@math.columbia.edu  \endemail
\thanks The first author was partially supported by NSF grant
DMS-92-03940. The second author was partially supported by NSF grant
DMS-94-02988.
\endthanks
\endtopmatter


\document

\section{Introduction}

Recently, Seiberg and Witten have introduced new invariants for smooth
$4$-manifolds which have led to dramatic progress in understanding the
$C^\infty$ properties of algebraic surfaces. Just as with Donaldson theory, the
new invariants are computed from a moduli space which, in case the underlying
$4$-manifold is a K\"ahler surface $X$, can be identified with a moduli space
of
holomorphic objects. In Donaldson theory, the holomorphic moduli space is the
space of holomorphic structures on a fixed $C^\infty$ complex vector bundle
over $X$ satsifying an additional nondegeneracy condition, stability. Such
moduli spaces have a rich geometric structure even for very simple K\"ahler
surfaces, such as $\Pee ^2$, and seem to become more progressively complicated
as
the surface becomes more complicated. In Seiberg-Witten theory, the relevant
moduli spaces are the spaces of complex curves $D$ on $X$, which are thus
parametrized by the Hilbert scheme of $X$, such that $D$ satisfies an
additional
numerical condition akin to stability. Now the structure of the Hilbert scheme
of curves on a smooth surface is an interesting problem in algebraic geometry.
However, it turns out for rather trivial reasons involving the Hodge index
theorem that the geometric interest of the Seiberg-Witten moduli spaces of a
surface $X$ is in a certain sense inversely proportional to the interest in $X$
itself as an abstract surface. Thus for example if $X$ is a minimal surface of
general type the Seiberg-Witten moduli spaces are two reduced points
corresponding to the trivial (empty) curve. Of course, it is this fact which
enables one to prove that the first Chern class of the canonical bundle of a
minimal surface of general type is a $C^\infty$ invariant up to sign. At the
other extreme, if $X$ is a ruled surface over a curve $C$ of genus at least
$2$,
then the Seiberg-Witten  moduli spaces are connected with the Brill-Noether
theory of special divisors on $C$, if $X= \Pee ^1\times C$ is a product ruled
surface, and to various interesting questions concerning stable bundles over
$C$
in general. Our goal in this paper is to discuss these and other related
examples.

The outline of this paper is as follows. In Section 1 we construct the Hilbert
scheme of a complex surface via $\dbar$ methods. To our knowledge, such a
construction has not appeared in the literature. In Section 2 we identify the
deformation complex for the Seiberg-Witten equations of a K\"ahler surface in
holomorphic terms and show that the Kuranishi model for the Seiberg-Witten
equations is the same as the Kuranishi model for the equations defining the
Hilbert scheme. In other words, the natural homeomorphism from the
Seiberg-Witten moduli space to the Hilbert scheme of ``stable" divisors on $X$
is
an isomorphism of real analytic spaces. In Section 3 we discuss how to make
computations in case the moduli space is smooth but does not have the expected
dimension, using the Euler class of the obstruction bundle. These arguments and
various generalizations are well-known to specialists in many different
contexts. In Section 4 we apply this study to elliptic surfaces. There is a
substantial overlap of the material in Sections 2--4 with the paper of
Brussee \cite{2}.

The remainder of the paper is concerned with ruled surfaces. We discuss the
infinitesimal and analytic structure of the moduli space for product ruled
surfaces in Section 5, and then  compute the invariant in
the special case where the curve involved is a section of the surface
(possibly with some fiber components). In Section 6, we deform the surface to
a general ruled surface and show that the Hilbert scheme of sections is much
better behaved: it is always smooth of the correct dimension. Using this
result, we give another computation of the invariants in the case of a section.
This computation goes back to Corrado Segre in 1889 \cite{11} and was given a
modern proof, for the case of the $0$-dimensional invariant, by Ghione
\cite{6}.
(Note that Segre considered the case of moduli spaces of
sections of arbitrary dimension.) We shall give a quick description of these
and related results. These methods generalize to compute the invariant in
general homologically; the problem is that it is not known whether, for a
general ruled surface, the Hilbert scheme always has the correct dimension.
Finally we remark that the computation of the invariant is a special case of
the
transition formula for Seiberg-Witten invariants for $4$-manifolds with $b_2^+
=
1$. This formula has been computed by the authors, by methods quite reminiscent
of those in Section 6, as well as by Li and Liu \cite{9}. Thus our goal in
Sections 5 and 6 has been, not so much to compute the invariant (although it is
amusing to see the connections with the enumerative calculations of
Brill-Noether theory) as it has been to see the relationship between the study
of
the Seiberg-Witten moduli spaces for ruled surfaces and questions in
Brill-Noether theory as well as the theory of rank two stable bundles on
curves.

\section{1. Structure of the Hilbert scheme.}

Let $X$ be an algebraic (or complex) surface, and let $D_0$ be an effective
divisor on $X$. We do not assume that $D_0$ is smooth or even reduced. Let
$H_{D_0,X}$ be  the Hilbert scheme of all effective divisors $D$ on
$X$ such that $c_1(\scrO_X(D))= c_1(\scrO_X(D_0))$ in $H^2(X; \Zee)$. As a set,
$H_{D_0,X}$ consists of all effective divisors $D$ homologous to $D_0$, i\.e\.
algebraically equivalent to $D_0$. There is a morphism
$H_{D_0,X}\to
\operatorname{Pic}X$ whose fibers are projective spaces. Over $X\times
H_{D_0,X}$ there is a tautological divisor $\Cal D$ whose restriction to each
slice $X\times \{t\}$ is the divisor $D_t$ on $X$ corresponding to $t$. A
general reference for the construction of $H_{D_0,X}$ and its properties is
\cite{10}.

The infinitesimal structure of $H_{D_0,X}$ is given as follows: from the
natural exact sequence
$$0 \to \scrO_X \to \scrO_X(D_0) \to \scrO_{D_0}(D_0)\to 0,$$
we have the associated long exact cohomology sequence. The Zariski tangent
space
to
$H_{D_0,X}$ is naturally the space of sections of the normal bundle
$H^0(D_0;\scrO_{D_0}(D_0))$. Note that the long exact cohomology sequence
gives
$$\gather
0 \to  H^0(X; \scrO_X(D_0))/H^0(X; \scrO_X) \to
H^0(D_0;\scrO_{D_0}(D_0))\\
\to H^1(X; \scrO_X)\to H^1(X; \scrO_X(D_0)).
\endgather$$
Here $H^0(X; \scrO_X(D_0))/H^0(X; \scrO_X)$ is the space of sections of
$\scrO_X(D_0)$ modulo the line through a nonzero section vanishing along $D_0$,
and is thus naturally the tangent space to the linear system $|D_0|$ at $D_0$.
The map $H^0(D_0;\scrO_{D_0}(D_0)) \to
H^1(X; \scrO_X)$ represents the infinitesimal change in the line bundle
$\scrO_X(D_0)$. We let $K_0$ denote the image of
$H^0(D_0;\scrO_{D_0}(D_0))$ in $H^1(X; \scrO_X)$, so that $K_0$ is the kernel
of the map from $H^1(X; \scrO_X)$ to $H^1(X; \scrO_X(D_0))$ defined by $\sigma
_0$. Thus there is an exact sequence
$$0 \to  H^0(X; \scrO_X(D_0))/H^0(X; \scrO_X) \to H^0(D_0;\scrO_{D_0}(D_0)) \to
K_0 \to 0.$$

The obstruction space to the deformation theory of $H_{D_0,X}$ is given as
follows: let $K_1$ be the image of $H^1(X; \scrO_X(D_0))$ in
$H^0(D_0;\scrO_{D_0}(D_0))$, or in other words the cokernel of the map from
$H^1(X; \scrO_X)$ to $H^1(X; \scrO_X(D_0))$ defined by $\sigma _0$. Then
$K_1$ is the obstruction space to the functor corresponding to $H_{D_0,X}$. If
$K_1=0$, then
$H_{D_0,X}$ is scheme-theoretically smooth at $D_0$ of dimension equal to $\dim
H^0(D_0;\scrO_{D_0}(D_0))$. (The converse is not necessarily true.) We say that
$D_0$ is {\sl semiregular\/} if $K_1=0$, or in other words if the map $H^1(X;
\scrO_X(D_0)) \to H^0(D_0;\scrO_{D_0}(D_0))$ is zero.

The following theorem was proved by Kodaira-Spencer
\cite{7} in the semiregular case (and was claimed by Severi):

\theorem{1.1} Let $D_0$ be a curve on $X$. Then the Zariski tangent space $T$
to
$H_{D_0,X}$ at $D_0$ fits into an exact sequence
$$0 \to  H^0(X; \scrO_X(D_0))/H^0(X; \scrO_X) \to T \to K_0 \to 0,$$
where $K_0 = \Ker \{\times \sigma _0\: H^1(X; \scrO_X)\to H^1(X;
\scrO_X(D_0))\}$.  Locally analytically in a neighborhood of $D_0$, $H_{D_0,X}$
is defined by the vanishing of a convergent power series without constant or
linear term from
$T$ to
$$K_1 =
\Coker \{\times \sigma _0\: H^1(X; \scrO_X)\to H^1(X; \scrO_X(D_0))\}.$$
\endstatement

To prove Theorem 1.1, one can analyze the deformation theory and obstruction
theory for $H_{D_0,X}$ via power series as in \cite{7} and
\cite{10}, and apply Schlessinger's theory. Here we give a $C^\infty$ proof of
Theorem 1.1. Given
$D_0$, let $L_0$ denote the $C^\infty$ complex line bundle defined by
$\scrO_X(D_0)$. From this point of view, the scheme $H_{D_0,X}$ is the set
(with
real analytic structure) of all $C^\infty$ sections of $L_0$ which are complex
analytic for some choice of holomorphic structure on $L_0$, modulo the action
of
the nowhere zero functions acting by multiplication. We fix a given holomorphic
structure on
$L_0$ with $\dbar$-operator simply denoted by $\dbar$, and a given nonzero
holomorphic section $\sigma _0$ of $L_0$ for this holomorphic structure.

The equations which say that $s$ is a holomorphic section
for some holomorphic structure on $L_0$ read as follows: there exist a
$\dbar$-closed $(0,1)$-form $A$ such that $(\dbar +A)(s) = 0$. Thus
$H_{D_0,X}$ is the zero set of the function
$$F_0\:  \Ker \dbar\oplus \Omega ^0(L_0)\subset \Omega
^{0,1}(X) \oplus \Omega ^0(L_0)\to  \Omega ^{0,1}(L_0)$$
defined by
$$F_0(A,s) = (\dbar +A)( s) = \dbar s +  As,$$
modulo the action of $\Cal G_\Cee$, the {\sl complex gauge group\/}, where
$\Cal
G_\Cee$ is the multiplicative group of nowhere vanishing $C^\infty$ functions
on
$X$ and $\lambda \in \Cal G_\Cee$ acts on
$(A,s)$ via $(A - \dbar \lambda/\lambda, \lambda \cdot s)$.
An easy calculation shows that $F_0\circ \lambda = \lambda
F_0$. Of course, in order to analyze the equations, we need to pass to an
appropriate Sobolev completion of all of these spaces, but we shall leave the
details of this standard procedure to the reader.

Next we calculate the linearized complex. The space
$\Omega ^0(L_0)$ of all
$C^\infty$ sections of $L_0$ is a vector space, and we may thus identify the
tangent space to $\Omega ^0(L_0)$ at a given section $\sigma _0$ vanishing at
$D_0$ with $\Omega ^0(L_0)$ again. The space of all $(0,1)$-conections on
$L_0$ is an affine space over $\Omega ^{0,1}(X)$, with origin the
$\dbar$-operator corresponding to the given complex structure, and so its
tangent space is $\Omega ^{0,1}(X)$. The nowhere zero functions on
$X$ may be (locally) identified with $\Omega ^0(X)$, the set of all $C^\infty$
functions on $X$, via the exponential, and the differential at $s= \sigma _0,
\lambda = 0$ of
$$\lambda \in \Omega ^0(X) \mapsto e^\lambda \cdot s$$
is multiplication by $\sigma _0$: $\lambda \mapsto \lambda \cdot \sigma _0$.
Taking the differential of the $\Cal G_\Cee$-action, we obtain a complex $\Cal
C_0$:
$$0 \to \Omega ^0(X) @>{d_1}>> \Ker \dbar\oplus\Omega ^0(L_0)  @>{d_2}>>
\Omega ^{0,1}(L_0).$$ Here the map $d_1\:\Omega ^0(X) \to \Ker
\dbar\oplus\Omega
^0(L_0)$ sends $\lambda$ to $(-\dbar \lambda,\lambda \cdot \sigma _0)$ and
$d_2\:\Ker \dbar\oplus\Omega ^0(L_0)  \to \Omega ^{0,1}(L_0)$ sends $(A,s)$ to
$\dbar s + A\sigma _0$. However, this complex is not elliptic. Thus, the
restriction of $F_1$ to $(\operatorname{Im} d_1)^\perp$ is not Fredholm.

To remedy the above problem, consider instead the function $F(A,s) = \pi \circ
F_0(A,s)$, where $\pi \: \Omega ^{0,1}(L_0) \to \Ker \dbar$ is
orthogonal projection onto the kernel of $\dbar$. (Note that $F$ is not in fact
equivariant with respect to the action of $\Cal G_\Cee$.)

\lemma{1.2} For $A$ in a neighborhood of zero in $\Ker \dbar$, $F(A,s) = 0$ if
and only if $F_0(A,s) = 0$.
\endstatement
\proof Clearly, if $F_0(A,s) =0$, then $F(A,s) =0$. Conversely, suppose that
$F(A,s) = 0$. This says that $F_0(A,s)$ is orthogonal to $\Ker \dbar$, so that
$F_0(A,s)= (\dbar +A)(s)=\dbar ^*\gamma$ for some $\gamma \in
\Omega ^{0,2}(L_0)$. By hypothesis $\dbar A =0$, so that $(\dbar +A)^2 =0$.
Applying $\dbar +A$ to $F_0$, we obtain $(\dbar +A)\dbar ^*\gamma =0$, or in
other words $\dbar \dbar ^*\gamma +A\dbar ^*\gamma = 0$. We claim that, in this
case, if $A$ lies in some neighborhood of zero, then $\dbar ^*\gamma =0$ and
thus $F_0(A,s) = 0$. We may restrict $\gamma $ to $(\Ker \dbar ^*)^\perp =
\operatorname{Im}\dbar$, and in this case $\dbar \dbar ^*$ is an isomorphism
from $\operatorname{Im}\dbar$ to itself (after taking appropriate completions).
Likewise, if $\pi_0$ denotes orthogonal projection from $\Omega ^{0,2}(L_0)$ to
$\operatorname{Im}\dbar$, then $\pi _0\circ( \dbar \dbar ^* + A\dbar ^*)$ is a
bounded map from
$\operatorname{Im}\dbar$ to itself, after taking appropriate completions, which
is invertible for $A=0$ and so for $A$ in a neighborhood of zero. It follows
that for $A$ in some neighborhood of zero, and for an arbitrary $\gamma \in
\Omega ^{0,2}(L_0)$, if $\dbar \dbar ^*\gamma +A\dbar ^*\gamma = 0$, then for
$\gamma _0$ the projection of $\gamma $ to $(\Ker \dbar ^*)^\perp$, we have
$\dbar ^*\gamma = \dbar ^*\gamma _0$ and $\pi _0\circ( \dbar \dbar ^* + A\dbar
^*)(\gamma _0) = 0$, so that $\gamma _0 = 0$ and $\gamma \in \Ker \dbar ^*$.
Hence
$\dbar ^*\gamma =0$ and so $F_0(A,s) =0$ as claimed.
\endproof

The linearization of the
equation $F$ and the gauge group action at
$(0, \sigma _0)$ gives a complex
$\Cal C_1$ defined by the top line of the following commutative diagram:
$$\CD
0@>>> \Omega ^0(X) @>{e_1}>> \Ker \dbar \oplus \Omega ^0(L_0)  @>{e_2}>> \Ker
\dbar @>>> 0\\ @. @|  @VVV  @VVV @.\\
@. \Omega ^0(X) @>>>\Omega ^{0,1}(X) \oplus \Omega ^0(L_0)  @>>> \Omega
^{0,1}(L_0), @.
\endCD$$
where the vertical maps are the natural inclusions, and the differentials
are given by
$e_1 (\lambda) = (-\dbar \lambda, \lambda \sigma _0)$ and
$e_2(A,s) = \dbar s  + A\sigma _0$. There is a subcomplex $\Cal C'$ defined by
$$ \Omega ^0(L_0) @>{\dbar}>> \Ker \{\,\dbar \: \Omega ^{0,1}(L_0)\to \Omega
^{0,2}(L_0)\,\},$$  shifted up a dimension, with differential $\dbar$, and the
quotient complex is isomorphic to the complex $\Cal C''$ defined by:
$$\Omega ^0(X)  @>{\dbar}>> \Ker \{\,\dbar \: \Omega ^{0,1}(X)\to \Omega
^{0,2}(X)\,\}.$$ Thus the deformation complex is elliptic and so the
restriction
of
$F$ to a slice for the $\Cal G_\Cee$ action is Fredholm. Taking the long exact
cohomology sequence associated to the exact sequence of complexes
$$ 0 \to \Cal C' \to \Cal C_1 \to \Cal C'' \to 0,$$
we see that the cohomology of
$\Cal C_1$ fits into the exact sequence
$$0 \to H^0(\scrO_X) \to H^0(L_0) \to H^1(\Cal C_1) \to H^1(\scrO_X) \to
H^1(L_0)
\to H^2(\Cal C_1) \to 0.$$
A routine calculation shows that the induced maps $H^i(\scrO_X) \to H^i(L_0)$
are given by multiplication by
$\sigma _0$. Thus $H^1(\Cal C_1)$ satisfies the exact sequence for $T$ given in
Theorem 1.1 and $H^2(\Cal C_1) \cong K_1$.  This concludes the proof of Theorem
1.1. \qed
\medskip

The following identifies the quadratic term of the obstruction map:

\proposition{1.3} Let $\xi \in H^0(\scrO_{D_0}(D_0))$ be an element of the
Zariski tangent space to $H_{D_0,X}$. Then the quadratic term in the
obstruction map is equal to $\delta \xi \cup \xi \in H^1(\scrO_{D_0}(D_0))$,
where $\delta \: H^0(\scrO_{D_0}(D_0)) \to H^1(\scrO_X)$ is the coboundary
map in the natural long exact sequence.
\endstatement
\proof It is easy to give the quadratic term of the obstruction map
using the power series approach of \cite{7}. In terms of the approach outlined
here, the quadratic term of $F(A,s)$ is $\pi (A\cdot s)$, where $\dbar A = 0$
and
$A\sigma _0 = -\dbar s$. If we identify the class of $(A,s)$ in $H^1(\Cal C_1)$
with an element $\xi \in H^0(\scrO_{D_0}(D_0))$, then it is easy to see that
the class of $A$ in $H^{0,1}(X) = H^1(\scrO_X)$ is exactly
$\delta \xi$. One then checks that the projection of $A\cdot s$ to
$\Ker \dbar$ corresponds to
$\delta \xi \cup \xi$.
\endproof

Note that, if we apply $\delta\: H^1(\scrO_{D_0}(D_0)) \to H^2(\scrO_X)$ to
the element $\delta \xi \cup \xi$, we obtain $\delta \xi \cup \delta \xi \in
H^2(\scrO_X)$, which is zero since $\delta \xi \cup \delta \xi = -\delta \xi
\cup \delta \xi$ as $\scrO_X$ is a sheaf of commutative rings. Thus $\delta \xi
\cup \xi$ lies in the image of $H^1(X; \scrO_X(D_0))$ in
$H^0(D_0;\scrO_{D_0}(D_0))$.

There is also clearly a universal divisor $\Cal D \subset X\times F^{-1}(0)$
defined by the vanishing of $s$. This completes the analytic construction of
$H_{D_0,X}$ and the discussion of semiregularity. Note that we have not
strictly speaking shown that $\Cal D$ is a divisor on the complex space
$X\times H_{D_0,X}$. This would need a discussion of relative $\dbar$-operators
similar to, but easier than, the discussion in \cite{4}, Chapter IV, 4.2.3. In
other words, we would need to show that $\Cal D$ is a Cartier divisor in the
possibly nonreduced complex space $X\times H_{D_0,X}$, which follows by showing
that locally on $X\times H_{D_0,X}$, there is a holomorphic embedding of the
complex space $X\times H_{D_0,X}$ in $X\times \Cee ^N$ for some $N$ so that
$\Cal D$ is locally the restriction of a complex hypersurface. Finally, to
identify this construction with the usual construction of
$H_{D_0,X}$, and to make a geometric identification of $H_{D_0,X}$ possible,
we would have to show that $H_{D_0,X}$ has a universal property. In other
words, given a complex space $T$, not necessarily reduced, and a Cartier
divisor on $X\times T$, flat over $T$, we need to exhibit a morphism of
complex spaces from $T$ to $H_{D_0,X}$. This again can be done along the lines
of  \cite{4}, Chapter IV. In the cases described in this paper, $H_{D_0,X}$
will
be smooth or a union of generically reduced components, and the arguments
needed
are substantially simpler than the arguments in the general case.

The divisor $\Cal D$ is a Cartier divisor and so there is a holomorphic line
bundle $\scrO_{X\times  H_{D_0,X}}(\Cal D)$ over $X\times  H_{D_0,X}$. Slant
product with $c_1(\scrO_{X\times  H_{D_0,X}}(\Cal D))$ defines a map $H_0(X)
\to H^2( H_{D_0,X})$, and we let $\mu$ be the image of the natural generator of
$H_0(X; \Zee)$ under this map. Clearly $\mu = \pi _2{}_*c_1(\scrO_{X\times
H_{D_0,X}}(\Cal D))= \pi _2{}_*[\Cal D]$. For a fixed $p\in X$, there is the
inclusion of the slice
$\{p\}\times  H_{D_0,X}$ in $X\times  H_{D_0,X}$, and clearly $\mu$ is the
first Chern class of the line bundle $\scrO_{X\times  H_{D_0,X}}(\Cal
D)|\{p\}\times  H_{D_0,X}$, under the natural identification of $\{p\}\times
H_{D_0,X}$ with $H_{D_0,X}$. If $\Cal D$ meets $\{p\}\times H_{D_0,X}$
properly, or in other words if there is no component $\Cal M$ of $ H_{D_0,X}$
such that $p$ lies in every divisor in $\Cal M$, then $\Cal D \cap \{p\}\times
H_{D_0,X}$ is a Cartier divisor in $\{p\}\times  H_{D_0,X} \cong H_{D_0,X}$,
whose support is the set of divisors $D$ such that $p\in D$, and this divisor
is a geometric representative for $\mu$. In fact, the divisor $\mu$ is an ample
divisor on $ H_{D_0,X}$, which can be shown  for example by using the method of
Chow schemes described in
\cite{10}, Lecture 16, and identifying the numerical equivalence class of $\mu$
up to a positive rational multiple with the natural ample divisor on the Chow
scheme.

There is another description of the complex line bundle corresponding to $\mu$.
For $p\in X$, let $\Cal G_\Cee ^0 \subset \Cal G_\Cee$, the {\sl based gauge
group\/}, be the set of
$\lambda \in \Cal G_\Cee$ such that $\lambda (p) =1$. Thus the quotient of
$\Cal G_\Cee$ by $\Cal G_\Cee^0$ is $\Cee ^*$, and if instead of dividing
out $F^{-1}(0)$ by the local action of $\Cal G_\Cee$ we divide out by $\Cal
G_\Cee^0$, the result is a $\Cee ^*$-bundle over $H_{D_0,X}$, which thus
corresponds to a complex line bundle $\Cal L_0(p)$. We claim that this line
bundle has first Chern class equal to $\mu$. First note that there is a
universal
$C^\infty$ complex line bundle $\Cal L_0$ over $X\times H_{D_0,X}$ whose
restriction to the slice $\{p\}\times H_{D_0,X}$ is $\Cal L_0(p)$. Here $\Cal
L_0$ is defined as follows: let $\Cal A_\Cee^*(L_0)$ be the set of pairs
$(A,s)$
where $A$ is a $(0,1)$-connection on $L_0$ and $s$ is a nonzero section of
$L_0$.  Then $\Cal G_\Cee$ acts freely on $\Cal A_\Cee^*(L_0)$; let the
quotient
be denoted $\Cal B_\Cee(L_0)$. Since $\Cal G_\Cee$ also acts as a group of
automorphisms of $L_0$, there is a line bundle $\Cal L_0$ over $X\times \Cal
B_\Cee(L_0)$ obtained by dividing out
$L_0 \times  \Cal A_\Cee^*(L_0)$ by the action of $\Cal G_\Cee$. We also denote
the restriction of this line bundle to $X\times F^{-1}(0)$ by $\Cal L_0$. So we
must identify $\Cal L_0$ with
$c_1(\scrO_{X\times  H_{D_0,X}}(\Cal D))$ (at least on the reduction of
$H_{D_0,X}$). The point is that the tautological section $(s, (A,s))$ of the
pullback of $L_0$ to $X\times \Cal A_\Cee^*(L_0)$ is $\Cal G_\Cee$-equivariant
and so descends to a section of $\Cal L_0$ over $X\times \Cal B_\Cee(L_0)$.
Restricting to $X\times F^{-1}(0)$, we see that $\Cal L_0$ has a section
vanishing at $\Cal D$, and this identifies $\Cal L_0$ with $\scrO_{X\times
H_{D_0,X}}(\Cal D)$.

\section{2. Deformation theory for Seiberg-Witten moduli spaces of K\"ahler
surfaces.}

In this section we recall the description of Seiberg-Witten moduli spaces for
K\"ahler surfaces and compare this description to the discussion of the Hilbert
scheme in the previous section. For general references on  Seiberg-Witten
moduli
spaces  of K\"ahler surfaces, see \cite{12}, \cite{3}, as well as
\cite{5}.

For a given metric $g$ on $X$ and Spin${}^c$ structure $\xi$ on $X$ with
determinant $L$, the (unperturbed) Seiberg-Witten equations for a pair $(A,
\psi)$, where $A$ is a connection on $L$ and $\psi$ is a section of $\Bbb
S^+(\xi)$, the plus spinor bundle associated to $\xi$, are
$$\align
\dirac_A\psi&=0;\\
F_A^+ =q(\psi) &=\psi\otimes
\psi^*-\frac{|\psi|^2}{2}\operatorname{Id}.
\endalign$$
We let $\Cal M_g(\xi)$ be the corresponding moduli space. In case $X$ is a
K\"ahler surface and $g$ is a K\"ahler metric with associated K\"ahler form
$\omega$, $\Bbb S^+(\xi) \cong \Omega ^0(L_0)
\oplus \Omega ^{0,2}(L_0)$ for a complex line bundle $L_0$, and $\Bbb S^-(\xi)
\cong \Omega ^{0,1}(L_0)$. Moreover $L_0 ^{\otimes 2} = L\otimes K_X$, where
$L = \det \xi$.  We assume that $\omega \cdot c_1(L) \neq 0$, so that there are
no reducible solutions, and for simplicity we fix $\omega \cdot c_1(L) <0$. In
this case, writing $\psi$ in components $(\alpha, \beta)$, where $\alpha$ is
a section of $\Omega ^0(L_0)$ and $\beta$ is a section of $\Omega
^{0,2}(L_0)$, the Seiberg-Witten equations become
$$\align
F^{0,2} &= \dbar A^{0,1}
=\bar\alpha\beta ;\\
(F_A^+)^{1,1} &=\frac{i}{2}(|\alpha|^2-|\beta|^2)\omega ;\\
\bar\partial_A\alpha+\bar\partial^*_A\beta &=0.
\endalign$$
Under the assumption that $\alpha \neq 0$, the equations $\dbar A^{0,1} =
\bar\alpha\beta$ and $\dbar _A(\bar\partial_A\alpha+\bar\partial^*_A\beta) = 0$
imply that $\beta = 0$, and that  $A$ is a
$(1,1)$-conection on $L$ \cite{12, 3, 5}. Hence $L$ and $L_0$ have given
holomorphic structures, and $\dbar _A\alpha = 0$, so that $\alpha$ is a
nonzero holomorphic section of $L_0$. Thus $\alpha$ defines an effective
divisor $D$ with $L_0 = \scrO_X(D)$. Taking gauge equivalence defines $\alpha$
up to scalars, or in other words as an element of $|D|$. Thus to each element
of $\Cal M_g(\xi)$, there is a well-defined element of $H_{D_0,X}$ for some
fixed divisor $D_0$ such that the $C^\infty$ line bundle underlying
$\scrO_X(2D_0) \otimes K_X^{-1}$ is $L$. Conversely, to every point of
$H_{D_0,X}$ we can associate an irreducible solution of the  Seiberg-Witten
equations mod gauge equivalence, in other words a point of $\Cal M_g(\xi)$,
which essentially follows from a theorem of Kazdan-Warner. It is easy to see
that
the map from
$\Cal M_g(\xi)$ to
$H_{D_0,X}$ is a homeomorphism, and we shall show that it is an isomorphism of
real analytic spaces in a suitable sense. As in the previous section, we shall
pass to Sobolev completions of all of the spaces of $C^\infty$ sections
involved
without making the choice of completions explicit.

We begin by discussing the deformation complex associated to the Seiberg-Witten
equations for a K\"ahler surface. For a general Riemannian 4-manifold $X$, at
an irreducible solution $(A_0, \psi)$ to the Seiberg-Witten
equations, the appropriate deformation complex $\Cal C$ is
$$0 \to i\Omega ^0(X;\Ar) @>{\delta _1}>> i\Omega ^1(X;\Ar) \oplus \Bbb
S^+(\xi)@>{\delta _2}>> i\Omega ^2_+(X; \Ar)
\oplus \Bbb S^-(\xi) \to 0.$$
Here $\Omega ^2_+(X)$ is the space of $C^\infty$ self-dual 2-forms. The
differentials are as follows:
$\delta _1(\lambda) = (-2d\lambda, \lambda \psi )$ and
$$\delta _2(A, \eta) = (d^+A - Dq_{\psi}(\eta), \dirac\eta + \frac12A\cdot
\psi).$$ Here $d^+$ is the self-dual part of $d$, $Dq_{\psi}$ is the
differential of the quadratic map $q$ in the SW equations, evaluated at $\psi$
on $\psi$, and
$\frac12A\cdot \psi$ is the linear term of $\dirac _{A+A_0}\psi $. A
calculation shows that
$$Dq_\psi(\eta) = \eta \otimes \psi ^* + \psi \otimes \eta ^* -
\operatorname{Re} \langle \psi, \eta \rangle \operatorname{Id}.$$
In general, it seems to be somewhat difficult to analyze
this complex. In the case of a K\"ahler surface $X$, however, we can give a
very explicit description of the cohomology of the deformation complex. First
we recall the notation of the previous section:
$$\align
K_0 &= \Ker \{\times \sigma _0\:
H^1(\scrO_X) \to H^1(L_0)\};\\
K_1 &= \operatorname{Coker}\{\times \sigma _0\: H^1(\scrO_X) \to
H^1(L_0)\};\\
K_2 &= \Ker \{\times \sigma _0\: H^2(\scrO_X) \to
H^2(L_0)\}.
\endalign$$

(Note that $\times \sigma _0\: H^2(\scrO_X) \to
H^2(L_0) $ is surjective.)

\theorem{2.1} Suppose that $X$ is a K\"ahler surface with K\"ahler metric $g$
and associated K\"ahler form $\omega$. Let $\sigma _0$ be a nonzero holomorphic
section of $L_0$ corresponding to an irreducible solution $(\sigma _0, 0)$ of
the Seiberg-Witten equations. Then the Zariski tangent space $H^1(\Cal C)$ to
the
Seiberg-Witten moduli space sits in an exact sequence
$$0 \to H^0(L_0)/\Cee \sigma _0 \to H^1(\Cal C) \to K_0\to 0$$ and the
obstruction space $H^2(\Cal C)$ sits in an exact sequence
$$0 \to K_1 \to
H^1(L_0)\} \to H^2(\Cal C) \to K_2\to 0.$$
\endstatement
\proof The complex $\Cal C$ has the following complex as its symbol
complex:
$$0 \to i\Omega ^0(X;\Ar) @>{(d,0)}>> i\Omega ^1(X;\Ar) \oplus \Bbb
S^+(\xi)@>{(d^+, \dirac)}>> i\Omega ^2_+(X; \Ar)
\oplus \Bbb S^-(\xi) \to 0.$$
Thus it is elliptic and its (real) index is the same as the index of the above
complex, namely
$$1-b_1(X) + b_2^+(X) - 2(h^0(L_0) + h^2(L_0) - h^1(L_0)) = 2\chi (\scrO_X) -
2\chi (L_0).$$
Here we have used the identification of $\dirac$ with $\dbar + \dbar ^*$ up to
a factor of $\sqrt{2}$. Note that
$\delta _1(\lambda ) = 0$ if and only if
$\lambda$ is constant and $\lambda \sigma _0 = 0$. Thus $H^0(\Cal C) = 0$,
which just says that the point $(A, (\sigma _0, 0)$ is an irreducible
solution to the Seiberg-Witten equations,  and we must identify the terms
$H^1(\Cal C)$ and
$H^2(\Cal C)$. Identify
$i\Omega^1(X; \Ar)$ with
$\Omega ^{0,1}(X)$, $\Bbb S^+(\xi)$ with $\Omega ^0(L_0) \oplus \Omega
^{0,2}(L_0)$,
$i\Omega _+^2(X; \Ar)$ with $i\Omega ^0(X; \Ar)\omega \oplus \Omega ^{0,2}(X)$,
and $\Bbb S^-(\xi)$ with $\Omega ^{0,1}(L_0)$. Under these identifications,
for
$\lambda \in i\Omega ^0(X; \Ar)$,
$\delta _1 \lambda = (-\dbar \lambda, \lambda \sigma _0, 0)$ which has image
inside $\Omega ^{0,1}(X)\oplus \Omega ^0(L_0)$. Moreover, $\delta _2 (A^{0,1},
\alpha, \beta)$ is given by
$$( (\dbar \bar A^{0,1} + \partial A^{0,1})^+ -
\operatorname{Re} \langle \sigma _0,  \alpha \rangle i\omega,
\dbar A^{0,1} -\bar \sigma _0 \beta,
\dbar \alpha + \dbar ^*\beta + \frac12A^{0,1} \sigma _0 ).$$ After identifying
$i\Omega ^0(X; \Ar)\omega$ with $\Omega ^0(X; \Ar)$ by taking $-i$ times the
contraction $\Lambda$ with $\omega$, we can write this as
$$\delta _2 (A^{0,1}, \alpha, \beta) = (T_1(A^{0,1}, \alpha), T_2(A^{0,1},
\alpha) + S(\beta)),$$ where (since $\omega \wedge \omega $ is twice the volume
form)
$$\align T_1(A^{0,1}, \alpha) &= -i(\Lambda\dbar \bar A^{0,1} + \Lambda
\partial A^{0,1}) - 2\operatorname{Re} \langle \sigma _0,  \alpha \rangle  \in
\Omega ^0(X; \Ar) ;\\ T_2(A^{0,1}, \alpha) &= (\dbar A^{0,1}, \dbar \alpha +
\frac12A^{0,1} \sigma _0) \in
\Omega ^{0,2}(X) \oplus \Omega ^{0,1}(L_0);\\ S(\beta) &= ( -\bar \sigma _0
\beta,\dbar ^*\beta) \in \Omega ^{0,2}(X) \oplus \Omega ^{0,1}(L_0).
\endalign$$ Thus $\delta _2(A^{0,1}, \alpha, \beta) =0$ if and only if
$T_1(A^{0,1}, \alpha) = 0$ and
$T_2(A^{0,1}, \alpha) + S(\beta) =0$. Next we claim:

\lemma{2.2} $T_2(A^{0,1}, \alpha) + S(\beta) =0$ if and only if $\beta = 0$ and
$T_2(A^{0,1}, \alpha) =0$. Moreover,  $\dbar A^{0,1} =\bar \sigma _0 \beta$ and
$\dbar (\dbar \alpha + \dbar ^*\beta + \frac12A^{0,1} \sigma _0) = 0$ if and
only if
$\beta = 0$ and $\dbar A^{0,1} =0$.
\endstatement
\proof Clearly, if $\beta = 0$ and
$T_2(A^{0,1}, \alpha) =0$, then $T_2(A^{0,1}, \alpha) + S(\beta) =0$.
Conversely, suppose that $T_2(A^{0,1}, \alpha) + S(\beta) =0$, or in other
words that $\dbar A^{0,1} =\bar \sigma _0 \beta$ and
$\dbar \alpha + \dbar ^*\beta + \frac12A^{0,1} \sigma _0=0$. Taking $\dbar$ of
the second equation, we find that
$$\dbar A^{0,1}\cdot\sigma _0 + \dbar \dbar ^*\beta = |\sigma _0|^2\beta +
\dbar
\dbar ^*\beta = 0.$$ Taking the inner product with $\beta$ shows that
$$\int _X|\sigma _0|^2 |\beta |^2 + \|\dbar ^*\beta\|^2 = 0.$$ Hence $\beta =
0$, and clearly then $T_2(A^{0,1}, \alpha) =0$ as well. The proof of the
second assertion is similar.
\endproof

Now we exhibit an isomorphism from $H^1(\Cal C)$ to $H^1(\Cal C_1)$, where
$\Cal C_1$ is the complex defined in the previous section, up to a factor of
$2$ (which arises because in our point of view $A^{0,1}$ is a connection on
$L_0^{\otimes 2} \otimes K_X$ rather than on $L_0$): $H^1(\Cal C_1)$ is the
quotient of
$$\{\, (A^{0,1}, \alpha): \dbar A^{0,1} = 0, \dbar \alpha +
\frac12A^{0,1}\sigma _0 = 0\,\}$$
by the subgroup of elements of the form $(-2\dbar f, f\sigma_0)$, where $f$ is
a {\sl complex\/} valued $C^\infty$ function. To exhibit this isomorphism,
given a class in $H^1(\Cal C)$ represented by $(A^{0,1}, \alpha, \beta)$, then,
by Lemma 2.2, $\beta = 0$ and $\dbar \alpha +
\frac12A^{0,1} \sigma _0 = 0$. Moreover $(A^{0,1}, \alpha)$ is well-defined up
to the subgroup of the form $(-2\dbar \lambda, \lambda\sigma_0)$, where
$\lambda$ is a purely imaginary $C^\infty$ function. Thus $(A^{0,1}, \alpha)$
defines an element of $H^1(\Cal C_1)$.

Conversely, start with  a representative $(A^{0,1}, \alpha)$ for an element of
$H^1(\Cal C_1)$. Then $(A^{0,1}, \alpha, 0)$ satisfies $T_2(A^{0,1},\alpha)
=0$ but not necessarily $T_1(A^{0,1}, \alpha) = 0$. On the other hand, we can
change $(A^{0,1}, \alpha)$ by an element of the form $(-2\dbar h, h\sigma
_0)$, where $h$ is a {\sl real\/} $C^\infty$ function, without affecting
$T_2(A^{0,1},\alpha)$. If we set $\gamma = T_1(A^{0,1}, \alpha)$, then
$$T_1(A^{0,1}-2\dbar h, \alpha + h\sigma _0) = \gamma - 2i(\Lambda \dbar
\partial h + \lambda \partial \dbar h) -2|\sigma _0|^2 h.$$
{}From the K\"ahler identities, $\Lambda \dbar = -i\partial ^*$ and similarly
$\Lambda \partial = i\dbar ^*$. Thus
$$i(\Lambda \dbar
\partial h + \lambda \partial \dbar h) = 2\operatorname{Re}\dbar ^*\dbar h =
-\Delta h,$$
where $\Delta$ is the negative definite Laplacian on $X$, and we seek to
solve the equation
$$2\Delta h -2|\sigma _0|^2 h  = -\gamma.$$
(Note that this equation is the linearized version of the Kazdan-Warner
equation used in identifying the Seiberg-Witten moduli space with the Hilbert
scheme.) Now we have the following:

\lemma{2.3} The operator $\Delta - |\sigma _0|^2$ is an isomorphism from
$\Omega ^0(X;\Ar)$ to itself.
\endstatement
\proof If $\Delta h- |\sigma _0|^2h =0$, then taking the inner product with
$h$ we find that $\|dh\|^2 = |\sigma _0|^2h ^2 =0$. Thus $h$ is constant and
$|\sigma _0|^2h ^2 =0$, so that $h=0$. Hence the operator
$\Delta - |\sigma _0|^2$ is injective. It is an elliptic operator on $\Omega
^0(X; \Ar)$ whose index is the same as the index of the Laplacian on functions,
namely zero. Thus it is also surjective.
\endproof

Thus given the initial representative $(A^{0,1},\alpha)$, there is a unique
choice of $h$ such that $T_1(A^{0,1}-2\dbar h, \alpha + h\sigma _0) = 0$.
Mapping $(A^{0,1},\alpha)$ to $(A^{0,1}-2\dbar h, \alpha + h\sigma _0, 0)$
then gives a well-defined map from $H^1(\Cal C_1)$ to $H^1(\Cal C)$, and
clearly the maps constructed are inverses. We have therefore showed that
$H^1(\Cal C) \cong H^1(\Cal C_1)$. By the proof of Theorem 1.1, there is an
exact sequence as claimed in the statement of Theorem 2.1.

We turn now to the identification of $H^2(\Cal C)$. Given $\gamma \in \Omega
^0(X;\Ar)$, it follows from Lemma 2.3 that we can solve the equation $\Delta h-
|\sigma _0|^2h = \gamma$. Thus there exists an $h$ such that $T_1(-2\dbar h,
h\sigma _0, 0) = \gamma$, and moreover $(-2\dbar h,h\sigma _0, 0)$ is in the
kernel of $T_2+S$. We can therefore identify the cokernel of $\delta _2$ with
the
cokernel of
$$T_2 + S\: \Omega ^{0,1}(X) \oplus \Omega ^0(L_0) \oplus \Omega ^{0,2}(L_0)
\to \Omega ^{0,2}(X) \oplus \Omega ^{0,1}(L_0).$$
Let $\Cal K$ denote the image of $T_2 + S$, so that $\Cal K$ is the set
$$\{\, (\dbar A^{0,1} - \bar \sigma _0\beta, \dbar \alpha +\dbar^*\beta +
\frac12A^{0,1}\sigma _0): A^{0,1} \in \Omega ^{0,1}(X) , \alpha \in \Omega
^0(L_0), \beta \in \Omega ^{0,2}(L_0)\,\}.$$
First consider the subgroup $0 \oplus \Ker \dbar \subset \Omega ^{0,2}(X)
\oplus \Omega ^{0,1}(L_0)$. For an element in $\Cal K\cap (0 \oplus \Ker
\dbar)$
which is the image of $(A^{0,1}, \alpha, \beta)$, we have
$\dbar A^{0,1} =\bar \sigma _0 \beta$ and
$\dbar (\dbar \alpha + \dbar ^*\beta + \frac12A^{0,1} \sigma _0) = 0$. By
Lemma 2.2, this condition is equivalent to
$\beta = 0$ and $\dbar A^{0,1} =0$. Thus
$$(0 \oplus \Ker \dbar)/\Cal K\cap (0 \oplus \Ker \dbar) \cong \Ker \dbar
/\{\,\dbar \alpha +\frac12A^{0,1} \sigma _0: \dbar A^{0,1} =0\,\}.$$
This is clearly the same as $H^1(L_0)/\sigma _0\cdot H^1(\scrO_X)= K_1$. So we
have found the subspace of $H^2(\Cal C)$ described in Theorem 2.1. The quotient
$K_2'$ of
$H^2(\Cal C)$ by $K_1$ is the same as $\Omega ^{0,2}(X) \oplus
\operatorname{Im}\dbar/(\operatorname{Id}\oplus \dbar)(\Cal K)$. Now
$(\operatorname{Id}\oplus \dbar)(\Cal K)$ is the subgroup
$$\{\, (\dbar A^{0,1} - \bar \sigma _0\beta, \dbar \dbar^*\beta +
\frac12\dbar A^{0,1}\sigma _0): A^{0,1} \in \Omega ^{0,1}(X) ,  \beta \in
\Omega ^{0,2}(L_0)\,\}.$$
If we consider the projection of this subgroup to the factor
$\operatorname{Im}\dbar \subseteq \Omega ^{0,2}(L_0)$, it is surjective since
$\dbar \dbar^*$ is an isomorphism on $\operatorname{Im}\dbar$. Thus $K_2'$ is
isomorphic to
$$\Omega ^{0,2}(X)/\{\, \dbar A^{0,1} - \bar \sigma _0\beta: \dbar \dbar^*\beta
+ \frac12\dbar A^{0,1}\sigma _0) = 0\,\}.$$
Let $\Cal K' = \{\, \dbar A^{0,1} - \bar \sigma _0\beta: \dbar \dbar^*\beta +
\frac12\dbar A^{0,1}\sigma _0) = 0\,\}$. We claim:

\lemma{2.4} $\Cal K' \cap \Ker \dbar ^* = \{\, -\bar \sigma _0\beta: \dbar
^*\beta = 0\,\}$.
\endstatement
\proof Suppose that
$$\dbar ^*\dbar A^{0,1} = \dbar ^*\bar \sigma _0\beta; \qquad
\dbar \dbar^*\beta =-
\frac12\dbar A^{0,1}\sigma _0.$$
Taking the inner product with $A^{0,1}$, we find:
$$\gather
\|\dbar A^{0,1}\|^2 = \langle \dbar ^*\dbar A^{0,1}, A^{0,1}\rangle
= \langle \dbar ^*\bar \sigma _0\beta, A^{0,1}\rangle = \langle \bar
\sigma _0\beta, \dbar A^{0,1}\rangle \\= \langle \beta, \sigma _0\dbar
A^{0,1}\rangle
= -2\langle \beta, \dbar \dbar^*\beta \rangle = -2\|\dbar^*\beta\|^2.
\endgather$$
It follows that $\dbar A^{0,1}=0$ and that $\dbar^*\beta =0$, and that  $\Cal
K' \cap \Ker \dbar ^*$ is as claimed.
\endproof

Using Lemma 2.4, there is an injection of $\Ker \dbar ^*/(\Cal K' \cap \Ker
\dbar
^*)$ into $K_2'$. Now $\Ker \dbar ^*\subseteq \Omega ^{0,2}(X)$ is
naturally $H^2(\scrO_X)$, and $\{\, -\bar \sigma _0\beta: \dbar ^*\beta =
0\,\} = \{\, -\sigma _0^*\beta: \dbar ^*\beta = 0\,\}$ is the image of
$H^2(L_0)$ under $\sigma _0^*$. The quotient $H^2(\scrO_X)/\sigma
_0^*H^2(L_0)$ is isomorphic to the orthogonal complement of
$\operatorname{Im}\sigma _0^*$, namely the kernel of multiplication by $\sigma
_0$ on $H^2(\scrO_X)$, which we have denoted by $K_2$. So there is an injection
of $K_2$ into $K_2'$. The real index of
$\Cal C$ is
$-\dim _\Ar H^1(\Cal C) + \dim _\Ar H^2(\Cal C)$, and we have shown that one
half the real index is at least
$$-h^0(L_0)+1 -\dim _\Cee K_0+ \dim _\Cee K_1+\dim_\Cee K_2.$$
Moreover equality holds only if $K_2'$ is exactly equal to $\Ker \dbar
^*/(\Cal K'
\cap \Ker \dbar ^*)$, and thus is isomorphic to
$K_2 =\Ker \{\times \sigma _0\:
H^2(\scrO_X) \to H^2(L_0)\}$.
On the other hand, the above alternating sum is the same as $\chi
(\scrO_X) - \chi (L_0)$ (take the alternating sum of the dimensions in the
cohomology exact sequence associated to multiplying by $\sigma _0$), which as
we have seen is one half the real index of $\Cal C$. It follows that $K_2' =
\Ker \dbar ^*/(\Cal K' \cap \Ker \dbar ^*)\cong K_2$ and that we have the
desired
exact sequence for $H^2(\Cal C)$.
\endproof

\corollary{2.5} If in the above notation $X$ is a K\"ahler surface, $L=
K_X^{-1}$ with the
\rom{Spin}${}^c$ structure corresponding to the trivial line bundle $L_0$, then
the Zariski tangent space is zero-dimensional and the obstruction space is
zero.
\endstatement
\proof In this case $H^0(L_0) = \Cee \sigma _0$, $\times \sigma _0\:
H^1(\scrO_X) \to H^1(L_0)$ is an isomorphism, and $\sigma _0H^0(K_X) =
H^0(K_X)$. Thus both the Zariski tangent space and the obstruction space are
zero.
\endproof

Next we compare the Kuranishi model of the Seiberg-Witten moduli space to the
Kuranishi model of the Hilbert scheme described in the previous section. In
what follows we assume that $H^2(\scrO_X) = 0$. In fact, if $X$ is a minimal
surface with $H^2(\scrO_X) \neq 0$, then either $X$ is of general type or it is
elliptic, a $K3$ surface, or a complex torus. In case $X$ is of general type,
the relevant Seiberg-Witten moduli spaces are smooth points corresponding to
$\pm K_X$ of the appropriate dimension, and the Kuranishi obstruction space is
zero by Corollary 2.5 above. In case $X$ is elliptic, the moduli space need
not be of the expected dimension, and the Kuranishi obstruction space need not
be zero, but we shall see in the next section that the obstruction map is
always identically zero and hence that the map  from $\Cal M_g(\xi)$ to
$H_{D_0,
X}$ is a diffeomorphism between two smooth manifolds. The other cases
involve reducible solutions to the Seiberg-Witten equations, and thus are
slightly exceptional from our point of view.  The case of a
nonminimal surface may then be reduced to the minimal case, at least for an
open set of K\"ahler metrics; we omit the details. Thus essentially the only
interesting case to consider is the case where
$H^2(\scrO_X) = 0$.

\theorem{2.6} Suppose that $H^2(\scrO_X) = 0$. Then the natural homeomorphism
from $\Cal M_g(\xi)$ to
$H_{D_0, X}$ is an isomorphism of real analytic spaces.
\endstatement
\proof  We keep the convention that $A$ is a connection on $L =
L_0^{\otimes 2}\otimes K_X$, rather than on $L_0$, and that it induces a
connection on $L_0$ once we have fixed once and for all a Hermitian connection
on
$K_X$. Recall that
$H_{D_0, X}$ is locally defined as the zeroes of the Fredholm map $F(A,\alpha)
=
\pi_{\Ker \dbar}(\dbar A + \frac12A\cdot \alpha)$, restricted to a slice of the
complex gauge group action on $\Ker \dbar \oplus \Omega ^0(L_0) \subset \Omega
^{0,1}(L_0) \oplus \Omega ^0(L_0)$.

As for $\Cal M_g(\xi)$, it is locallly defined by the zero set of the three
equations $G = (F_A^+)^{1,1} -\frac{i}{2}(|\alpha|^2-|\beta|^2)\omega$, $\dbar
A^{0,1}-\bar\alpha\beta$, and
$\bar\partial_A\alpha+\bar\partial^*_A\beta$. Setting the first equation
$G$ equal to zero on a slice $S'$ for the real gauge group gives a slice for
the
complex gauge group: indeed, the differential of $G$ is the map $T_1$ defined
in the proof of Theorem 2.1, and Lemma 2.1 shows that, given $(A^{0,1},
\alpha)$, there is a unique real-valued $C^\infty$ function $h$ such
that $T_1(A^{0,1}-2\dbar h, \alpha + h\sigma _0) =0$.  Thus
$$T_1^{-1}(0) \oplus iT\Cal G_\Ar = \Omega
^{0,1}(X) \oplus \Omega ^0(L_0) \oplus \Omega ^{0,2}(L_0).$$
In particular $S' \cap G^{-1}(0)$ is a slice for the complex gauge group in a
neighborhood of the origin; denote this slice by $S$.

Consider now the remaining two equations. Defining
$$\Cal F(A, \alpha, \beta) = (\dbar A^{0,1}-\bar\alpha\beta,
\bar\partial_A\alpha+\bar\partial^*_A\beta),$$
we can view $\Cal F$ as a section of the trivial vector bundle over $\Omega
^{0,1}(X) \oplus \Omega ^0(L_0) \oplus \Omega ^{0,2}(L_0)$ with fiber $\Omega
^{0,2}(X) \oplus \Omega ^{0,1}(L_0)$ whose restriction to the slice $S$ is
Fredholm and locally defines $\Cal M_g(\xi)$.

By our assumption that $H^2(\scrO_X) = 0$, and since  $H^2(\scrO_X)$ surjects
onto $H^2(L_0)$, it follows that $\dbar\: \Omega ^{0,1}(L_0) \to \Omega
^{0,2}(L_0)$ is surjective. Hence the natural map
$$\Omega ^{0,2}(X) \oplus \Omega ^{0,1}(L_0) \to \Ker \dbar \oplus \Omega
^{0,2}(X) \oplus \Omega ^{0,2}(L_0)$$
defined by $(\psi, \eta ) \mapsto (\pi _{\Ker \dbar}\eta, \psi, \dbar \eta)$,
is an isomorphism. Thus for small $A$, the map
$$(\psi, \eta ) \mapsto (\pi _{\Ker \dbar}\eta, \psi, \dbar _A\eta)$$ is again
an isomorphism. We may then view this map as an automorphism of the trivial
vector bundle over $\Omega
^{0,1}(X) \oplus \Omega ^0(L_0) \oplus \Omega ^{0,2}(L_0)$ (or an
appropriate neighborhood of the origin) with fiber
$\Omega ^{0,2}(X) \oplus \Omega ^{0,1}(L_0)$. Under this automorphism, $\Cal F$
corresponds to the section $(F_1, F_2)$, where
$$\align
F_1 &= \pi _{\Ker \dbar}(\dbar _A\alpha + \dbar _A^*\beta);\\
F_2 &= ( \dbar A^{0,1}-\bar\alpha\beta,\dbar _A^2\alpha + \dbar _A\dbar
_A^*\beta).
\endalign$$
Thus the Kuranishi model for $\Cal F$ on the slice $S$ is the same as that for
the pair $(F_1, F_2)$ on $S$. It is easy to check that the differential of the
map $F_2$ is the same as the differential  of $\Cal F$ followed by $(\Id,
\dbar)$. In other words, the cokernel of the differential of $F_2$ is exactly
the group $K_2$ of Theorem 1.1, namely the kernel of multiplication from
$H^2(\scrO_X)$ to $H^2(L_0)$. Since $H^2(\scrO_X) = 0$ by assumption, the
differential of $F_2$ is onto and $F_2^{-1}(0)$ is a smooth submanifold of a
neighborhood of the origin in $\Omega
^{0,1}(X) \oplus \Omega ^0(L_0) \oplus \Omega ^{0,2}(L_0)$. Now $F_2 = 0$ if
and only if $\dbar A^{0,1}=\bar\alpha\beta$ and $\dbar _A^2\alpha + \dbar
_A\dbar
_A^*\beta = 0$. As we mentioned earlier, these equations imply that $\beta = 0$
and hence that $\dbar A^{0,1}=0$. Conversely, if $\beta = 0$ and $\dbar A^{0,1}
=
0$, then $F_2(A, \alpha, \beta ) = 0$. Solving the equation $\Cal F = 0$ on the
slice $S$ is the same then as solving the equation $\pi _{\Ker \dbar}(\dbar
_A\alpha) = 0$ on the slice $S\cap (\Ker \dbar \oplus \Omega ^0(L_0))$, at
least
in a neighborhood of the origin. This is exactly the equation $F(A, \alpha)$ on
the slice $S\cap (\Ker \dbar \oplus \Omega ^0(L_0))$ for the complex gauge
group. A standard argument (see for example \cite{4}, Chapter 4, proof of
Theorem 3.8) shows that this is the same as the usual Kuranishi model for $F$,
in other words that this model is isomorphic to the Kuranishi model formed by
taking any other slice for the $\Cal G_\Cee$-action. Thus the two Kuranishi
models are isomorphic as complex spaces.
\endproof

\section{3. Obstruction bundles.}

Fix an oriented  $4$-manifold $X$ with a Riemannian metric $g$ (not necessarily
a
K\"ahler surface). Let
$\Cal A^*(L)$ denote the spaces of pairs $(A, \psi)$, where
$A$ is a connection on $L$ and $\psi$ is a nonzero section of $\Bbb S^+(\xi)$
as
in the previous section. The real gauge group $\Cal G$ acts on $\Cal A^*(L)$,
and we denote the quotient by $\Cal B(L)$. The trivial Hilbert space bundle
$i\Omega ^2_+(X; \Ar) \oplus \Bbb S^-(\xi)$ descends to a Hilbert bundle $\Cal
H$ over
$\Cal B(L)$, and the moduli space $\Cal M_g(\xi)$ is the zero set of the
Fredholm
section $F(A, \psi)$ defined by the Seiberg-Witten equations. As such $\Cal
M_g(\xi)$ has a real analytic structure and in particular a Zariski tangent
space. In this section, we are concerned with the following situation:
suppose that the space $\Cal M_g(\xi)$ is a smooth compact manifold, not
necessarily of the expected dimension. Thus the dimension of the Zariski
tangent
space of
$\Cal M_g(\xi)$ at every point is equal to the dimension of $\Cal M_g(\xi)$ at
that point, and these tangent spaces fit together to form the tangent bundle
$T\Cal M_g(\xi)$ to
$\Cal M_g(\xi)$. Note that the tangent bundle is in fact just $\Ker dF\: T\Cal
A^*(L) \to \Cal H$. It follows that the obstruction spaces have locally
constant
rank on $\Cal M_g(\xi)$ and thus, by standard elliptic theory, fit together to
form a vector bundle $\Cal O$ over $\Cal M_g(\xi)$. In case $g$ is a K\"ahler
metric on the complex surface $X$, the arguments of the previous section show
that the fiber of $\Cal O$ over a point $(A_0, \sigma _0)$ may be canonically
identified with the middle cohomology of the elliptic complex
$$\Omega ^{0,1} (X; \Ar) \oplus \Omega ^0(L_0) \to \Omega ^{0,2} (X; \Ar)
\oplus
\Omega ^{0,1}(L_0) \to  \Omega ^{0,2}(L_0),$$
where the first map is $(A^{0,1} , \alpha) \mapsto (\dbar A^{0,1}, \dbar \alpha
+ \frac12A^{0,1}\sigma _0)$ and the second map is $(\varphi, \psi) \mapsto
\dbar \psi - \frac12\varphi \cdot \sigma _0)$. Again by a slight modification
of the standard theory for the $\dbar$-operator, it follows that the bundle
$\Cal O$ is a holomorphic vector bundle  over $\Cal M_g(\xi)$. In fact,
letting $\frak C$ be the complex $\scrO_{X\times \Cal M_g(\xi)} \to \scrO_{X
\times \Cal M_g(\xi)}(\Cal D_0)$, it is easy to see that $\Cal O$ is
the $C^\infty$ vector bundle associated to the hyperdirect image $\Ar
^2\pi _2{}_*\frak C \cong R^1\pi _2{}_*\scrO_{\Cal D_0}(\Cal D_0)$.  To
use this information to evaluate Seiberg-Witten
invariants, we have the following:

\theorem{3.1} In the above notation, suppose that the expected real dimension
of
$\Cal M_g(\xi)$ is $2d$, and let $\mu$ be the natural class in $H^2(\Cal
B(L))$.
Then the value of the Seiberg-Witten function on $\xi$ is $\int _{\Cal
M_g(\xi)}
e(\Cal O)\cup
\mu ^d$, where $e(\Cal O)$ is the Euler class of the vector bundle $\Cal O$. In
case
$X$ is a K\"ahler surface and
$\Cal M_g(\xi)$ is equidimensional, this is the same as
$\int _{\Cal M_g(\xi)}c_n(\Cal O)\cup \mu ^d$, where $n =
\operatorname{rank}\Cal O$.
\endstatement

In particular we need to calculate $c_n(\Cal O)$:

\lemma{3.2} In the above notation, suppose that $X$ is a K\"ahler surface and
that $\Cal M_g(\xi)$ is smooth. Then $c_n(\Cal O)$ is the degree $n$ term in
$$c(\pi _2{}_!\scrO_{X \times \Cal M_g(\xi)}(\Cal D_0))^{-1}c(T\Cal
M_g(\xi)).$$
\endstatement
\proof We need to calculate $c_n(\Cal O) = c_n (R^1\pi _2{}_*\scrO_{\Cal
D_0}(\Cal D_0))$. Now the morphism $\Cal D_0 \to \Cal M_g(\xi)$ has relative
dimension one, and so $\pi _2{}_!\scrO_{\Cal D_0}(\Cal D_0)$, which is by
definition the alternating sum of the
$R^i\pi _2{}_*\scrO_{\Cal D_0}(\Cal D_0)$, is just $R^0\pi
_2{}_*\scrO_{\Cal D_0}(\Cal D_0) - R^1\pi _2{}_*\scrO_{\Cal D_0}(\Cal D_0)$.
Furthermore, $R^0\pi _2{}_*\scrO_{\Cal D_0}(\Cal D_0)$ is just the tangent
bundle $T\Cal M_g(\xi)$ to $\Cal M_g(\xi)$. Thus
$$c(\Cal O)= c(R^1\pi _2{}_*\scrO_{\Cal D_0}(\Cal D_0)) = c(\pi
_2{}_!\scrO_{\Cal D_0}(\Cal D_0))^{-1}c(T\Cal M_g(\xi)).$$
Using the exact sequence
$$0 \to \scrO_{X \times \Cal M_g(\xi)} \to \scrO_{X \times \Cal M_g(\xi)}(\Cal
D_0) \to \scrO_{\Cal D_0}(\Cal D_0) \to 0,$$
it follows that, in the $K$-theory of $\Cal M_g(\xi)$,
$$\pi _2{}_!\scrO_{\Cal D_0}(\Cal D_0) = \pi _2{}_!\scrO_{X \times \Cal
M_g(\xi)}(\Cal D_0) - \pi _2{}_!\scrO_{X \times \Cal
M_g(\xi)}.$$
Since $\pi _2{}_!\scrO_{X \times \Cal
M_g(\xi)}$ is a trivial vector bundle,
$$c(\pi _2{}_!\scrO_{\Cal D_0}(\Cal D_0)) = c(\pi _2{}_!\scrO_{X \times \Cal
M_g(\xi)}(\Cal D_0)).$$
Putting this together with the above formula for $c(\Cal O)$ gives the
statement of (3.2).
\endproof

\demo{Proof of \rom{(3.1)}} Consider quite generally the following situation:
$\Cal H \to \Cal B$ is a Hilbert vector bundle over the connected Hilbert
manifold $\Cal B$, and
$\sigma$ is a smooth section of $\Cal H$. Let $Z = \sigma ^{-1}(0)$, assumed
connected for the sake of simplicity, and suppose that the differential
$d\sigma$ is Fredholm of index
$e$, at least in a neighborhood of $Z$. Suppose moreover that $Z$ is a
smooth compact submanifold of $\Cal B$ of finite dimension $e'$ and that $\Ker
d\sigma _z$ has constant rank for every
$z\in Z$ and that the corresponding subbundle of
$T\Cal B|Z$ is the tangent bundle to $Z$. Define the obstruction bundle $\Cal
O\to Z$ by $\Cal O = \Coker (d\sigma |Z)$, of rank $e' - e$. Theorem 3.1 is
then a consequence of the following lemma, which implies that the class of a
small generic perturbation of $Z$ is the same as the class of a generic section
of $\Cal O$, in other words that the Euler class of $\Cal O$ represents the
same
cohomology class as the Seiberg-Witten class of a generic moduli space.
\enddemo

\lemma{3.3} In the above situation, suppose that $\sigma _1$ is a
small nonlinear Fredholm $C^1$ perturbation of $\sigma$, with $\sigma
_1^{-1}(0) = Z_1$, and that
$\sigma _1$ is transverse to $0$ in the sense that $d\sigma _1$ is surjective
at
every point of $Z_1$. Then there exists a section $s$ of $\Cal O \to Z$ which
is transverse to $0$ and a small isotopy in $\Cal B$ from $s^{-1}(0)$ to $Z_1$.
\endstatement
\proof Since $Z$ is compact, standard arguments show that there is a
neighborhood $\nu$ of $Z$ in $\Cal B$ which is diffeomorphic to a Hilbert disk
bundle over $Z$. Let $\pi \: \nu \to Z$ be the retraction. Over $Z$, there is
an orthogonal splitting $\Cal H|Z \cong \operatorname{Im}d\sigma \oplus \Cal
O$. Using $\pi$, we can pull this decomposition back to a splitting of $\Cal
H|\nu \cong I \oplus \Cal O$. Let $\pi _1 \: \Cal H|\nu \to I$ and $\pi _2 \:
\Cal H|\nu \to \Cal O$ be the projections. Consider the composed map $\pi
_1\circ \sigma \: \nu \to I$. At $z\in Z$, the differential of this map is
just $d\sigma$, and so restricted to a fiber $\pi ^{-1}(z)$, the differential
is an isomorphism. It follows that, if
$\nu$ is sufficiently small, then $\pi_1\circ \sigma |\pi ^{-1}(z)$ is an open
embedding for all $z\in Z$.

Now let $\sigma _1$ be a small perturbation of $\sigma$. If $\sigma _1$ is
sufficiently close to $\sigma$, then $\sigma _1^{-1}(0) = Z_1 \subseteq \nu$.
Consider the map $\pi _1 \circ \sigma _1 \: \nu \to I$. If we restrict this map
to a fiber $\pi ^{-1}(z)$ of the map $\pi \: \nu \to Z$,  $\pi_1\circ \sigma _1
|\pi ^{-1}(z)$ is close to an open embedding. Thus, as long as $\sigma _1$ is
close to $\sigma$, $\pi_1\circ \sigma _1
|\pi ^{-1}(z)$ is also an open embedding. In particular, $(\pi_1\circ \sigma
_1)^{-1}(0) = \hat Z_1$ is again a section of $\pi \: \nu \to Z$, and it is
close to the zero section $\sigma$. Thus $\hat Z_1$ is isotopic to $Z$ via a
small isotopy in $\nu \subseteq \Cal B$.

Clearly $Z_1 = \sigma _1^{-1}(0) = (\pi _2 \circ \sigma _1|\hat
Z_1)^{-1}(0)$. Moreover $(\pi _2 \circ \sigma _1|\hat
Z_1)$ is a section $\hat s$ of the restriction $\Cal O \to \hat Z_1$, and, if
$\sigma _1$ is transverse to $0$, then $\hat s$ is also transverse to $0$.
Using the isotopy constructed above to identify the section $\hat s$ with a
section of $\Cal O\to Z$, we see that we have indeed identified $Z_1$ up to
isotopy with a transverse section of $\Cal O \to Z$, as claimed.
\endproof

There are obvious generalizations of Lemma 3.3, and so of Theorem 3.1. For
example, we might only assume that $Z$ is a stratified space with $K$ a compact
subset contained in an open subset $U$ of $Z$, such that $U$ is a smooth
manifold and $\Ker d\sigma _u = TU_u$ for all $u\in U$. Then a generic small
perturbation
$\sigma _1$ of
$\sigma$ has the property that there exists a neighborhood $N$ of $K$ such that
$\sigma _1^{-1}(0) \cap N$ is a smooth manifold isotopic to a transverse
section of the obstruction bundle over $N$. For example, we might take for $K$
a subset of the form $\mu _1 \cap \dots \cap \mu _k$, where the $\mu _i$ are
generic geometric representatives for the $\mu$-divisor. However, we shall not
try to formulate the most general possible result along these lines.

\section{4. Elliptic surfaces.}

Let $X$ denote an elliptic
surface. Suppose that $f$ is the divisor class of a general fiber, the
multiple
fibers are
$F_i$, and that the multiplicity of $F_i$ is $m_i$. Thus $K_X = (p_g-1)f +
\sum _i(m_i-1)F_i$. We first consider the much simpler case where $X$ is
regular,
since this case arises in the smooth classification of elliptic surfaces (see
for
example \cite{5}, \cite{2}). Then we will discuss the multiplicities for a
general elliptic surface.

If
$X$ is regular,
$H^1(\scrO_X) = 0$ and the Seiberg-Witten obstruction space involves only the
two terms $H^1(\scrO_X(D_0))$  and $H^2(\scrO_X(D_0))$. The divisor $D_0$ is
semiregular if and only if $H^1(\scrO_X(D_0)) = 0$. It follows from (2.3) of
\cite{5} that the $D_0$ are exactly the effective divisors  which are
numerically equivalent to $\frac{1-r}2K_X$, for a rational number $r \leq 1$.
In particular $D_0 \cdot K_X = 0$. Another way to describe the $D_0$ is that
they are the effective divisors numerically equivalent to a rational multiple
of $K_X$  such that $K_X - 2D_0$ has positive fiber degree (is a positive
rational multiple of the fiber, or equivalently of $K_X$). A similar statement
holds if $X$ is not necessarily assumed to be regular. As
$L^2 =0$, the expected dimension of the moduli space is always zero, i\.e\. $X$
is of simple type. We now compute the dimensions of the cohomology groups:

\lemma{4.1} Suppose that $X$ is regular, and let $D_0 = af + \sum _ib_iF_i$
with
$a\geq 0$ and
$0\leq b_i\leq m_i-1$. Then:
$$\align
h^0(D_0) &= a+1;\\
h^1(D_0) &= \cases 0, & \text{if $a\leq p_g$};\\
a-p_g, &\text{if $a\geq p_g$}.
\endcases\\
h^2(D_0) &= \cases p_g-a, & \text{if $a\leq
p_g$};\\ 0, &\text{if $a\geq p_g$}.
\endcases
\endalign$$
\endstatement
\proof The statements about $h^0(D_0)$ and $h^2(D_0)$ are clear and the rest
follows from Riemann-Roch, since $\chi (\scrO_X(D_0)) = 1+p_g$.
\endproof

Thus we see that $D_0$ is semiregular if and only if $a\leq p_g$. However,
since $X$ is regular and thus $H_{D_0, X}$ is equal
to the linear system $|D_0|$, which is a projective space $\Bbb P^a$, $H_{D_0,
X}$ is always smooth and the Zariski tangent space is the actual tangent space.
To calculate the value of the Seiberg-Witten invariant on $L$, we take the top
Chern class of the obstruction bundle. Now the moduli space is $X\times \Pee
^a$, where $\Pee ^a = |D_0|$. Over $X\times \Pee ^a$, there is the incidence
divisor $\Cal D$. The obstruction bundle over $\Pee ^a$ has two terms $R^1\pi
_2{}_*\scrO_{X\times \Pee ^a}(\Cal D)$ and $R^2\pi _2{}_* \scrO_{X\times \Pee
^a}(\Cal D))$.

\proposition{4.2} Suppose in the above notation that $a\leq p_g -1$. Then the
multiplicity of the Seiberg-Witten invariant is $(-1)^a\dsize
\binom{p_g-1}{a}$.
\endstatement
\proof By Lemma 4.1, the first term of the obstruction bundle is zero, and we
must compute the top Chern class $c_a$ of $R^2\pi _2{}_* \scrO_{X\times \Pee
^a}(\Cal D))$.
Next we calculate the class of $\Cal D$
in $\Pic (X\times \Pee ^a)\cong \pi _1^*\Pic X \oplus \Zee \cdot \pi
_2^*\scrO_{\Pee ^a}(1)$. Since $\Cal D$ is the incidence divisor, its
restriction to the slice $\pi _2^*\{D\}$ is the divisor $D$, whereas its
restriction to any slice $\{p\}\times \Pee ^a$ such that $p$ is not in the base
locus of $D_0$ is a hyperplane in $\Pee ^a$. Thus $\scrO_{X\times \Pee ^a}(\Cal
D) \cong \pi _1^*\scrO_X(D_0) \otimes \pi _2^*\scrO_{\Pee ^a}(1)$. So we have
$$R^2\pi _2{}_* \scrO_{X\times \Pee ^a}(\Cal D)) = R^2\pi _2{}_*\left(\pi
_1^*\scrO_X(D_0) \otimes \pi _2^*\scrO_{\Pee ^a}(1) \right)= \scrO_{\Pee
^a}(1)^{p_g -a}.$$
Setting $h = c_1(\scrO_{\Pee ^a}(1))$, we want to take the term of degree $a$
in
$$\left((1+h)^{p_g-a}\right)^{-1} = (1+h)^{-(p_g-a)}.$$
By the binomial theorem (see below for our conventions on binomial
coefficients), this is
$\dsize \binom{-(p_g -a)}{a} =
(-1)^a\binom{p_g-1}{a}$.
\endproof

A similar argument shows:

\proposition{4.3} Suppose in the above notation that $a\geq p_g$. Then the
multiplicity of the Seiberg-Witten invariant is $1$ if $p_g =0$ and is
otherwise $0$.
\endstatement
\proof Since $H^2( D) =0$, by Lemma 4.1, we seek
$c_a\left(R^1\pi _2{}_*\scrO_{X\times \Pee ^a}(\Cal D)\right)$.
Using the
calculation
$\scrO_{X\times \Pee ^a}(\Cal D) \cong \pi _1^*\scrO_X(D_0) \otimes \pi
_2^*\scrO_{\Pee ^a}(1)$ given above, we need to find $c_a$ of
$$R^1\pi _2{}_*\left(\pi _1^*\scrO_X(D_0) \otimes \pi
_2^*\scrO_{\Pee ^a}(1)\right) =H^1(D_0) \otimes _\Cee \scrO_{\Pee
^a}(1) = (\scrO_{\Pee ^a}(1))^{a-p_g}.$$ If $p_g =0$, then $c((\scrO_{\Pee
^a}(1))^a) = (1+h)^a$, and thus
$c_a((\scrO_{\Pee ^a}(1))^a) = 1$. Otherwise, the multiplicity is $c_a$ of a
bundle of rank less than $a$, and so is zero.
\endproof

\noindent {\bf Remark.} If $p_g =0$, then the multiplicity is always $1$,
although we can have $a>0$ if there are more than two multiple fibers. If $p_g
>0$ and there are at most two multiple fibers (the case of finite cyclic
fundamental group), then it is easy to check that $a\leq p_g-1$. In general
however both of the terms in the exact sequence for the obstruction bundle can
be nonzero.
\medskip

Next we turn to elliptic surfaces which are not necessarily regular.  To state
the result, let us record the following convention on binomial coefficients
(made so that the binomial theorem holds): For $m\geq 0$, the binomial
coefficient
$\dsize\binom{n}{m} =\frac{1}{m!}n(n-1) \cdots (n-m+1)$. Thus it is $1$ if
$m=0$, by the usual conventions on the empty product, and it is $0$ if
$0\leq n<m$ and $m\neq 0$. Moreover
$\dsize\binom{-n}m = (-1)^m\binom{n+m-1}{m}$.
For example, $\dsize\binom{-1}{m} = (-1)^m$. With this
said, suppose that
$X$ is a minimal elliptic surface over a smooth curve of genus $g$ and $D_0 =
af
+ \sum _ib_iF_i$ is a basic class for $X$. Then the value of the Seiberg-Witten
invariant on
$D_0$ is given by the following formula, due independently to Brussee
\cite{2}:

\proposition{4.4} Let $\pi \: X \to C$ be a minimal elliptic surface over a
smooth curve $C$ of genus $g$ and let $D_0$ be an effective divisor
corresponding to a Seiberg-Witten basic class. Suppose that $D_0 = \pi ^*\bold
d + \sum _i a_iF_i$, where $\bold d$ is an effective divisor on $C$ of degree
$d$, the $F_i$ are the multiple fibers on $X$, of multiplicity $m_i$, and $0
\leq a_i \leq m_i -1$. Then the multiplicity of the Seiberg-Witten invariant is
$\dsize (-1)^d\binom{\chi(\scrO_X) + 2g-2}{d}$.
\endstatement
\proof If $D_0 = \pi ^*\bold
d + \sum _i a_iF_i$ as above, there is a natural morphism from $\Sym ^dC = C_d$
to $H_{D_0, X}$ obtained by pulling back the universal divisor on $C\times C_d$
to $X\times C_d$. Slightly tedious arguments left to the reader show that this
identifies $H_{D_0, X}$ as a set with $C_d$. To calculate Zariski tangent
spaces, it is easy to see that $H^0(\scrO_{D_0}(D_0))$ has dimension $d$ by
using the fact that the normal bundle of $F_i$ is torsion of order exactly
$m_i$. Thus the dimension of the Zariski tangent space to
$H_{D_0, X}$ is equal to the dimension of $H_{D_0, X}=C_d$, and the map $C_d
\to H_{D_0, X}$ is an isomorphism.

Let $\pi _2 \: X \times C_d \to C_d$ be the second projection, let $p\: X\times
C_d \to C_d$ be the map induced by $\pi$, and let $p_2\: C \times C_d \to C_d$
be second projection. Over $C\times C_d$ we have the incidence divisor $\Cal I$
defined by $\Cal I = \{\, (t, \bold d): t\in \operatorname{Supp}\bold d\,\}$.
Thus the universal divisor $\Cal D$ on $X \times C_d$ is just $p^*\Cal I$. Let
$\psi _2 \: C\times C_d \to C_d$ be projection onto the second factor. By
flat base change, $R^i\pi _2 {}_*p^*\scrO_{C\times C_d}(\Cal I) = R^i \psi
_2{}_*\scrO_{C\times C_d}(\Cal I)$, and in particular this is zero for $i=2$.
Recall that the multiplicity of the Seiberg-Witten invariant is then given by
evaluating
$c_a(R^1\pi _2{}_*\scrO_{\Cal D}(\Cal D))$ over $C_d$. From the exact sequence
$$0 \to \scrO_{X\times C_d} \to  \scrO_{X\times C_d}(\Cal D) \to \scrO_{\Cal
D}(\Cal D) \to 0,$$
we see that, in the $K$-theory of $C_d$,
$$\pi _2{}_!\scrO_{\Cal D}(\Cal D) =\sum _i(-1)^iR^i\pi _2{}_*\scrO_{\Cal
D}(\Cal D)= R^0\pi _2{}_*\scrO_{\Cal D}(\Cal D) - R^1\pi _2{}_*\scrO_{\Cal
D}(\Cal D)$$ agrees with $\pi _2{}_!\scrO_{X\times C_d}(\Cal D)$ up to the
trivial element $\pi _2{}_!\scrO_{X\times C_d}$, and thus they have the same
Chern classes. Moreover
$$c(R^1\pi _2{}_*\scrO_{\Cal D}(\Cal D)) = c(\pi _2{}_!\scrO_{X\times C_d}(\Cal
D))^{-1}c(R^0\pi _2{}_*\scrO_{\Cal D}(\Cal D)).$$
Finally $R^0\pi _2{}_*\scrO_{\Cal D}(\Cal D)$ is just the tangent bundle
$T_{C_d}$ to
$C_d$. By \cite{1}, p\. 322,
$$c(T_{C_d}) = (1+x)^{d+1-g}e^{-\theta /1+x},$$
where $x$ is the class of the divisor $C_{d-1}\subset C_d$ and $\theta$ is the
pullback of the theta divisor on $\Pic^d C$ under the natural map. We also have
the formula
$$\theta ^kx^{d-k} = \frac{g!}{(g-k)!}.$$
To find $c(\pi _2{}_!\scrO_{X\times C_d}(\Cal
D))$, we first apply the Grothendieck-Riemann-Roch theorem to find $\ch (\pi
_2{}_!\scrO_{X\times C_d}(\Cal D))$:
$$\ch (\pi _2{}_!\scrO_{X\times C_d}(\Cal D)) = \pi _2{}_*\left(\ch
(\scrO_{X\times C_d}(\Cal D))\pi _1^*\Todd X\right).$$
If $\delta = [\Cal I]$ is the class of $\Cal I$ on $C\times C_d$, then $\ch
(\scrO_{X\times C_d}(\Cal D)) = p^*e^\delta$. Moreover $\Todd X = 1 + rf + \chi
(\scrO_X)\cdot \text{pt}$ for some rational number $r$. By \cite{1}, p\.
338, $\delta =n[\text{pt}]\otimes 1 + \delta ^{1,1} + 1\otimes x$, where $x$ is
the class defined above and  $(\delta ^{1,1})^2 = -2[\text{pt}]\otimes \theta$,
$(\delta ^{1,1})^2 = (\delta ^{1,1})\cdot [\text{pt}]\otimes 1 = 0$. Since $\pi
^*[\text{pt}] = f$ with $f^2 = 0$, an easy calculation shows that
$$\ch (\pi _2{}_!\scrO_{X\times C_d}(\Cal D))= \pi _2{}_*\left(\ch
(\scrO_{X\times C_d}(\Cal D))\pi _1^*\Todd X\right) = \chi(\scrO_X) \cdot e^x$$
and thus (setting $\chi(\scrO_X) = \chi$ for brevity)
$$c(\pi _2{}_!\scrO_{X\times C_d}(\Cal D)) = (1+x)^\chi.$$
Finally, then, the multiplicity of the Seiberg-Witten invariant is the term of
degree $d$ in
$$(1+x)^{-\chi}(1+x)^{d+1-g}e^{-\theta /1+x} =
(1+x)^{d+1-g-\chi}e^{-\theta /1+x}.$$
Let $N = d+1-g-\chi$. Then
$$\align
(1+x)^Ne^{-\theta /1+x} &= \sum _i\binom{N}{i}\sum _j\frac{1}{j!}(-1)^j\theta
^j\sum _k\binom{-j}{k}x^{i+k}\\
&= \sum _a\sum _{i+j+k= a}
\binom{N}{i}\frac{1}{j!}(-1)^j\binom{-j}{k}\frac{g!}{(g-j)!}\\
&=\sum _a\sum _{i+j+k= a} (-1)^{j+k}\binom{N}{i}\binom{j+k-1}{k}\binom{g}{j}.
\endalign$$
Thus the degree $d$ term is
$$\sum _j(-1)^j\left(\sum
_k(-1)^k\binom{N}{d-j-k}\binom{j+k-1}{k}\right)\binom{g}{j}.$$
To evaluate the term
in parentheses above, we have the straightforward combinatorial lemma:

\lemma{4.5} We have:
$$\sum _{k=0}^a \binom{a+j+e}{a-k}\binom{-j}{k} = \sum _{k=0}^a
(-1)^k\binom{a+j+e}{a-k}\binom{j+k-1}{k} =
\binom{a+e}{a}.$$
\rom(By our conventions on
binomial coefficients, this is
$1$ if $a=0$ and is zero for $-a \leq e<0$ and $a\neq 0$.\rom)
\endstatement
\proof This follows by comparing the coefficient of $t^a$ in the two
different power series expansions of $(1+t)^{a+e} = (1+t)^{a+e+j}(1+t)^{-j}$.
\endproof

Returning to the proof of (4.4), the lemma shows that the term in parentheses
is
$\dsize \binom{N-j}{N-j-e}$, where $e = 1-g-\chi$ (take $a=d-j$ and $N= a + j +
e$). Thus we obtain
$$\align
\sum _j(-1)^j\binom{N-j}{N-j-e}\binom{g}{j} &= \sum _j(-1)^j\binom{d-j+ 1
-g-\chi}{d-j}\binom{g}{j}
\\ &= (-1)^d\sum _j\binom{\chi + g -2}{d-j}\binom{g}{j},
\endalign$$
which is just $(-1)^d$ times the
coefficient of $t^d$ in  $(1+t)^{\chi + g -2}(1+t)^g=(1+t)^{\chi + 2g -2}$,
namely  $\dsize (-1)^d
\binom{\chi + 2g -2}{d}$, as claimed.
\endproof

\section{5. Product ruled surfaces.}

In this section we shall consider the ruled surfaces $X$ of the form $\Pee
^1\times C$, where $C$ is a curve of genus $g\geq 1$. We shall also always
assume that $C$ is a generic curve in the sense of Brill-Noether theory, and
shall use \cite{1} as a general reference for the theory of special
divisors on curves.

Let $\pi _1\: X = \Pee ^1\times C \to \Pee ^1$ be the projection onto the first
factor and let $\pi _2\: \Pee ^1\times C \to C$ be the projection onto the
second. Let $F_1 = \pi _2^{-1}\{p\}$ be a fiber isomorphic to $\Pee ^1$ and let
$F_2 = \pi _1^{-1}\{p\}$ be a fiber isomorphic to $C$. Thus $F_i^2 = 0$ and
$F_1\cdot F_2 = 1$. In general we shall refer to a divisor numerically
equivalent
to $nF_1 + mF_2$ as a divisor {\sl of type\/} $(n,m)$, and similarly for a
complex line bundle. Thus for example
$K_X$ is of type $(2g-2, -2)$. Hence $K_X^2 = -8(g-1)$. Let $L$ be a line
bundle of type $(2a, 2b)$, so that $c_1(L) \equiv K_X \mod 2$. Then $L^2 =
8ab$, and so $L^2 \geq K_X^2$ if and only if $ab \geq 1-g$. Next suppose that
$L_0 = (K_X\otimes L)^{1/2}$ has a holomorphic section for some holomorphic
structure on $L_0$. As $L_0$ is of type $(g-1+a, b-1)$, we must have $b\geq 1$
and $a\geq 1-g$, and we can write $L_0 = \scrO_X(D_0)$, where $D_0$ is
linearly equivalent to $(b-1)\pi _1^*(\text{pt}) + \pi _2^*\bold d$ for some
divisor
$\bold d$ on
$C$ of degree
$d = g-1 + a$.

Next we turn to the condition that $L\cdot \omega < 0$ for some K\"ahler form
$\omega$. The real cohomology classes of K\"ahler metrics are exactly the
classes $\omega $ of type $(x,y)$ with $x, y\in \Ar$, and $x,y>0$. Thus
$$\omega \cdot L = 2xb + 2ay.$$
Since $b\geq 1$, we must have $a<0$, and it is clear that by choosing $x/y <
-a/b$, we can then arrange $\omega \cdot L < 0$. (Note conversely that if
$x/y > -a/b$, then $L$ does not correspond to a basic class. Since $-a \leq
g-1$
and $b\geq 1$, if we choose $x/y \geq g-1$, then there are no basic classes.)
The final conditions on
$a$ and $b$ are:
$$b \geq 1; \qquad   \frac{1-g}b \leq a <0.$$
We note that the expected (complex) dimension of the Seiberg-Witten moduli
space
is $g-1 +ab$. However, as we shall see, the actual dimension is equal to the
expected dimension only for $b=1$.

Given a curve $D_0$ of type $(g-1 + a, b-1) = (d,b-1)$, its irreducible
components correspond to curves of type $(e_i, c_i)$ with $\sum _ie_i = d$ and
$\sum _ic_i = b-1$. For example, if $D_0$ is irreducible, then $D_0$ is
simultaneously a cover of $\Pee ^1$ of degree $d$ and a cover of $C$ of degree
$b-1$. If $b=1$, then necessarily $D_0$ is a union of $d$ copies
of $C$, or more precisely a divisor of the form $\pi _2^*\bold d$ for some
divisor $\bold d$ of degree $d$ on $C$. In this case,
$H_{D_0, X}$ is just $C_d$, the
$d^{\text{th}}$ symmetric product of $C$ with itself. Note that $-K_X$ is a
divisor of this type, with $a = 1-g$. In general, for such divisors, we have:

\proposition{5.1} In case $c_1(L) = 2aF_1 + 2F_2$, then $H_{D_0, X}= C_d$ and
the value of the Seiberg-Witten invariant is $1$.
\endstatement
\proof We have seen that $H_{D_0, X}= C_d$ as sets. There is an obvious
universal divisor on $X\times C_d$ which is the pullback of the universal
divisor on $C\times C_d$. Thus there is a morphism from $C_d$ to $H_{D_0, X}$.
To see that this morphism is an isomorphism of schemes, it will suffice to show
that $H_{D_0, X}$ is smooth of dimension $d$. It is an easy exercise to
identify $H^i(\scrO_{D_0}(D_0))$ with $H^i(\Pee ^1; H^0(\scrO_{\bold d})\otimes
\scrO_{\Pee ^1}) =  H^i(\Pee ^1;\scrO_{\Pee ^1}) \otimes H^0(\scrO_{\bold
d})$. This has dimension $d$ for $i=0$ and is zero for $i=1$. Thus $H_{D_0, X}$
is smooth of dimension $d$, and is therefore isomorphic to $C_d$. Note that $X$
is not of simple type if $d>0$.

Clearly, for $p\in X$ with $\pi _2(p) = t\in C$, the incidence divisor for
$H_{D_0, X}$ and $p\in X$ may be identified with the incidence divisor for
$C_d$ and $t\in C$, along with multiplicities. Let $x$ be the class of the
divisor in $C_d$. By choosing $d$ distinct points $t_1, \dots, t_d$ and
checking that the intersections are transverse, we see that $x^d = 1$. Thus
$\mu ^d = 1$ for the Seiberg-Witten moduli space as well, so that the value of
the invariant is 1.
\endproof

For the remainder of this section, we shall mainly be interested in the case
$b=2$. In this case every curve
$D_0$ of type $(d,1)$ can be written either as $D_0 = D_1 + \pi _2^*\bold d_2$,
where
$D_1$ is the graph of a map $C\to \Pee ^1$ of degree $d_1$ and $\bold d_2$ is a
divisor of degree $d_2 = d-d_1$ on $C$, or
$D_0 = \pi _1^*(\text{pt}) + \pi _2^*\bold d$, where
$\bold d$ is a divisor of degree $d$ on $C$. More generally, an irreducible
divisor
$D_0 \subset \Pee ^1 \times C$ of type $(d, b-1)$ which is a section of
the line bundle $\pi _1^*\scrO_{\Pee ^1}(b-1) \otimes \pi _2^*\scrO_C(\bold
d)$ corresponds to a map $C \to
\Sym^{b-1}\Pee ^1 \cong \Pee ^{b-1}$. In this case it is easy to check that
the pullback of $\scrO_{\Pee ^{b-1}}(1)$ to $C$ is just $\scrO_C(\bold d)$.
Moreover, let $V$ be the smallest linear subspace of $\Cee ^b$ such that $\Pee
(V) \subseteq \Pee^{b-1}$ contains the image of $C$. Then $V$ is naturally a
quotient of $H^0(\bold d)^*$, of dimension $r+1$, say, corresponding to a
linear subseries of $|\bold d|$.  Note in particular that we always have $r\leq
b-1$.

Next suppose that $D_0$ is not necessarily irreducible. Then $D_0$ still
corresponds to a linear system $V \subseteq H^0(\bold d)^*$ with $\Pee (V)
\subseteq \Sym^{b-1}\Pee ^1 \cong \Pee ^{b-1}$. In fact, if $D_0$ is defined by
the section $\sigma _0 \in H^0(\bold b) \otimes H^0(\bold d)$, write $\sigma _0
= \sum _i\alpha _i \otimes \beta _i$, where the $\alpha _i \in H^0(\bold
b)$ and the $\beta _i
\in H^0(\bold d)$ are linearly independent. For $p\in C$, the morphism $C\to
\Sym^{b-1}\Pee ^1$ sends
$p$ to $\sum _i\beta _i(p)\alpha _i$, after choosing a coordinate for
$\scrO_C(\bold d)$ at $p$. This is well-defined if $p$ is not in the base locus
of the span of the $\beta _i$, and extends to a unique morphism $C\to
\Sym^{b-1}\Pee ^1$. Now $\sigma _0\in  H^0(\bold b) \otimes H^0(\bold d)$
defines a homomorphism $H^0(\bold d)^* \to H^0(\bold b)$ whose image $V$ is
spanned by the $\alpha _i$. Thus $V^*\subseteq H^0(\bold d)$ is a linear
series.

At one extreme, consider divisors $D_0$ of the form $\pi _1^*\bold b + \pi
_2^*\bold d$, where $\deg \bold b= b-1$ and $\deg \bold d = d$. In this case,
$V$ has dimension 1, the linear series corresponding to $V$ is the single
divisor $\bold d$, which is the base locus for the series, and the map $C \to
\Sym^{b-1}\Pee ^1
\cong
\Pee ^{b-1}$ is constant, with image equal to $\bold b$. In this  case, let
$\frak M_0 = \Pee ^{b-1}\times C_d$, thinking of this space as
parametrizing all divisors of the form $\pi _1^*\bold b + \pi _2^*\bold d$.
Then
we have an obvious divisor on the product
$X\times \frak M_0$, and thus there is an injective morphism $\Pee ^{b-1}
\times
C_d\to H_{D_0, X}$. In particular, $H_{D_0, X}$ has dimension at least $d+b-1
= g-1 + a + b-1 = g+a+b-2$, whereas the expected dimension of $H_{D_0, X}$ is
$g+ab-1$. In
this case the difference between $\dim \frak M_0$ and the expected dimension of
$H_{D_0, X}$ is $(b-1)(1-a)$. Since $a<0$ we see that the actual dimension is
always greater than the expected dimension as long as $b>1$. To see the image
of
the tangent space to $\frak M_0$ inside the Zariski tangent space of $H_{D_0,
X}$, note that  the Zariski tangent space of $H_{D_0, X}$ is
$H^0(\scrO_{D_0}(D_0))$. In case $D_0 =  \pi _1^*\bold b + \pi _2^*\bold d$,
there is an map $H^0(\scrO_{\pi _1^*\bold b}) \oplus H^0(\scrO_{\pi
_2^*\bold d}) \to H^0(\scrO_{D_0}(D_0))$ given as follows: If we set $E_1 = \pi
_1^*\bold b$ and $E_2 = \pi _2^*\bold d$, then this is just the natural map
$$(\scrO_X(E_1)/\scrO_X) \oplus (\scrO_X(E_2)/\scrO_X) \to
\scrO_X(E_1+E_2)/\scrO_X.$$ Locally for
$R = \Cee\{z_1, z_2\}$, this is the same as the map
$$R/z_1^{a_1}R \oplus R/z_2^{a_2}R \to R/z_1^{a_1}z_2^{a_2}R$$
defined by $(f,g) \mapsto z_2^{a_2}f+ z_1^{a_1}g$, which is an inclusion since
$z_1^{a_1}$ and $z_2^{a_2}$ are relatively prime. Thus
the image of the tangent space of
$\frak M_0$ is
$L_0 = H^0(\scrO_{\pi _1^*\bold b}) \oplus H^0(\scrO_{\pi _2^*\bold d})$, of
dimension
$b-1 + d$, and the map from $\frak M_0$ to $H_{D_0, X}$ is an immersion.

Concerning the structure of $H_{D_0, X}$ at a divisor $D_0$ of the form $\pi
_1^*\bold b + \pi _2^*\bold d$, where $\deg \bold b = b-1$ and $\deg \bold d =
d = g-1 +a$, we have the following result:

\lemma{5.2} The scheme $H_{D_0, X}$ is smooth of dimension $(b-1)(1-a)$ at the
divisor $D_0 = \pi _1^*\bold b + \pi _2^*\bold d$ if and only if $h^0(C;
\bold d) = 1$, or in other words if and only if $\bold d$ is an effective
divisor which does not move in a nontrivial linear system. More generally, if
$\dim |\bold d| = r$ and if
$T$ is the Zariski tangent space of $H_{D_0, X}$ at $D_0$, then there is an
exact
sequence
$$ 0 \to (H^0(\Pee ^1; \bold b) \otimes H^0(C; \bold d))/\Cee \sigma _0 \to T
\to \Ker \{H^1(\scrO_C) \to H^1(\scrO_C(\bold d))\to 0,$$
where the map $H^1(\scrO_C) \to H^1(\scrO_C(\bold d))$ is given by
multiplying by the section coresponding to $\bold d$. Thus
$$\dim T = b(r+1) -1 + d-r.$$
Finally, the image of
$\frak M_0$ is a component of $H_{D_0, X}$.
\endstatement
\proof First note that $H^0(X; \scrO_X(D_0)) = H^0(\Pee ^1; \scrO_{\Pee
^1}(b-1)) \otimes H^0(C; \bold d)$. Thus if $h^0(C; \bold d) =1$, then every
nonzero section of $H^0(X; \scrO_X(D_0))$ is of the form $\pi _1^*\bold b + \pi
_2^*\bold d$ for some divisor $\bold b$ of degree $b-1$ on $\Pee ^1$. Hence
$h^0(D_0) = b$ and $|D_0|\cong \Pee ^{b-1}$. In general, setting $r+1 =
h^0(\bold d)$,
$$\dim H^0(\scrO_X(D_0))/\Cee\cdot \sigma _0 = b(r+1) - 1.$$
Now $H^1(\scrO_X) \cong H^0(\scrO_{\Pee ^1}) \otimes H^1(\scrO_C)$. Given a
section $\sigma _0 $ of $\scrO_X(D_0)$ of the form $\sigma _1\otimes \sigma
_2$,
multiplication by
$\sigma_0$ is just multiplication by $\sigma _2$ from $H^1(\scrO_C)$ to
$H^1(\bold d)$, followed by  multiplication by $\sigma _1$. On the other
hand, using the exact sequence
$$0 \to \scrO_C \to \scrO_C(\bold d) \to \scrO_{\bold d} \to 0,$$
where the map  $\scrO_C \to \scrO_C(\bold d)$ is multiplication by $\sigma _2$,
it follows that the dimension of the kernel of the map $H^1(\scrO_C)\to
H^1(\bold d)$ is
$d - r$ (and the dimension of the cokernel is $g-d + r$). So the dimension of
the
Zariski tangent space is $b(r+1) - 1 + d-r$, whereas the actual dimension of
$\frak M_0$ is $b-1+d$.  For the generic divisor $\bold d \in C_d$,
$r=0$, since $d = g+a-1 <g$. For such a divisor $\bold d$, the map $\frak
M_0 \to H_{D_0, X}$ is an embedding near $\bold d$. Thus the image of
$\frak M_0$ is a component of
$H_{D_0, X}$. Finally  $H_{D_0, X}$ is smooth at $\bold d \in \frak M_0$ if and
only if $r=0$.
\endproof

The same argument identifies the obstruction space at $D_0$:

\lemma{5.3} Suppose as above that $D_0 = \pi _1^*\bold b + \pi _2^*\bold d$.
Then $$\Coker \{\times \sigma _0\: H^1(\scrO_X) \to H^1(\scrO_X(D_0))\}$$ is
equal to $H^0(\scrO_{\Pee ^1}(b-1))/\Cee \cdot \sigma _1 \otimes H^1(\bold d)$
and has dimension $(b-1)(g-d + r)$. \qed
\endstatement

We can put Lemma 5.3 in  a more intrinsic global form as follows. The tangent
bundle to $\Pee ^{b-1}$ is naturally $\scrO_{\Pee ^{b-1}}(1)^b/\scrO_{\Pee
^{b-1}}$. Over $C\times C_d$ we have the incidence divisor $\Cal I$ defined by
$$\Cal I = \{\, (t, \bold d): t\in \operatorname{Supp}\bold d\,\}.$$
Let $\psi
_2 \: C\times C_d \to C_d$ be projection onto the second factor. Then using the
exact sequence
$$0 \to \scrO_{C\times C_d} \to \scrO_{C\times C_d}(\Cal I) \to \scrO_{\Cal
I}(\Cal I) \to 0,$$
we have an exact sequence of direct image sheaves by applying $R^i\psi _2{}_*$.
Here $\psi _2|\Cal I$ is a $d$-sheeted cover, and in particular it is finite.
Moreover $\Cal I$ is a hypersurface in $C\times C_d$ and $\Cal I \cap \psi
_2^{-1}(\bold d)$ is identified with the divisor $\bold d$ on $C$. It is easy
to see that $R^0\psi _2{}_*\scrO_{\Cal I}(\Cal I) $ is canonically the tangent
bundle $T_{C_d}$ of $C_d$. Since it is torsion free, and  $R^0\psi
_2{}_*\scrO_{C\times C_d} \to R^0\psi _2{}_*\scrO_{C\times C_d}(\Cal I)$ is an
isomorphism at a general $\bold d$, $R^0\psi _2{}_*\scrO_{C\times
C_d} \cong R^0\psi _2{}_*\scrO_{C\times C_d}(\Cal I)$ and we have an exact
sequence
$$0 \to T_{C_d} \to R^1\psi _2{}_*\scrO_{C\times C_d} \to R^1\psi
_2{}_*\scrO_{C\times C_d}(\Cal I) \to 0.$$
Here $R^1\psi _2{}_*\scrO_{C\times C_d} = H^1(\scrO_C) \times \scrO_{C_d}$ is a
trivial bundle of rank $g$ on $C_d$, and we set $\Cal E =  R^1\psi
_2{}_*\scrO_{C\times C_d}(\Cal I)$. Note that the Chern classes of $\Cal
E$ are given by $c(\Cal E) = c(T_{C_d})^{-1}$.

\lemma{5.4} Let $p_1\: \Pee ^{b-1} \times
C_d \to \Pee ^{b-1}$ be projection onto the first factor and let $p_2\: \Pee
^{b-1} \times C_d \to C_d$ be projection onto the second factor. We then
have a map $p_1^*T_{\Pee ^{b-1}} \otimes p_2^*\Cal E \to R^1\pi
_2{}_*\scrO_{X\times \frak M_0}(\Cal D)/R^1\pi _2{}_*\scrO_{X\times \frak
M_0}$,
and it is an isomorphism over those points $(\bold b, \bold d)$ of $\frak M_0$
where $h^0(\bold d) = 1$. \qed
\endstatement

Let us study the obstruction space for a divisor
$D_0$ which is not necessarily in
$\frak M_0$.

\lemma{5.5} Let $D_0 \subset \Pee ^1 \times C$ be a divisor of type $(d, b-1)$,
corresponding to a morphism $C \to \Sym^{b-1}\Pee ^1 \cong \Pee ^{b-1}$. Let
$V \subset H^0(\bold d)^*$ be the linear subspace corresponding to the image
of $C$. Then the obstruction space for
$H_{D_0, X}$  at $D_0$ is zero if and only if the map $C \to \Pee
^{b-1}$ is nondegenerate and the map
$$\mu _0 \: V\otimes H^0(K_C - \bold d) \to H^0(K_C)$$
given by cup product is injective, which holds for a generic curve $C$.
\endstatement
\proof The obstruction space is given as the cokernel of the map $H^1(\scrO_X)
\to H^1(\scrO_X(D_0))$ which is multiplication by $\sigma _0$. Here
$H^1(\scrO_X) \cong H^0(\scrO_{\Pee ^1}) \otimes H^1(\scrO_C)$ and
$H^1(\scrO_X(D_0)) \cong H^0(\scrO_{\Pee ^1}(b-1)) \otimes H^1(\scrO_C(\bold
d))$. Let $V$ be the quotient of $H^0(\bold d)^*$
corresponding to the map $\sigma _0 \in H^0(\bold b) \otimes H^0(\bold d) =
\Hom (H^0(\bold d)^*, H^0(\bold b))$, and let
$\alpha _1, \dots, \alpha _n$ be a basis for $V$ viewed as a subspace of
$H^0(\bold b)$. It follows that we can write
$\sigma _0 = \sum _i\alpha _i \otimes \beta _i$ for $\beta _i \in H^0(\bold
d)$ which are linearly independent. In this case, $\{\beta _i\}$ must also be a
basis for $V^*
\subseteq H^0(\bold d)$. Multiplication by $\sigma_0$ is
equivalent to the map sending $\xi \in H^1(\scrO_C)$ to $\sum _i\alpha
_i\otimes
(\beta _i\xi)$. This map is surjective if and only if $V = H^0(\scrO_{\Pee
^1}(b-1))$ and the natural map
$$H^1(\scrO_C) \to H^0(\scrO_{\Pee ^1}(b-1)) \otimes H^1(\bold d)$$
is surjective.  After applying an invertible element of $H^0(\scrO_{\Pee
^1}(b-1))$, we can assume that $\alpha _i = \beta _i^*$, the dual basis to
$\beta _i$. In this case multiplication by $\sigma _0$ is easily seen to be the
adjoint of the map $\mu _0$. Thus  multiplication by $\sigma _0$ is surjective
if and only if $V = H^0(\scrO_{\Pee ^1}(b-1))$ and $\mu _0$ is injective.
\endproof

We note that, in case $\mu _0$ is injective, and in particular for a
generic curve $C$, (5.5) identifies the obstruction space as $(H^0(\scrO_{\Pee
^1}(b-1))/V)\otimes H^1(\bold d)$. For example, in case $\sigma _0 = \pi
_1^*\bold b+ \pi _2^*\bold d$, $V$ is a line in
$H^0(\scrO_{\Pee ^1}(b-1))$ and the obstruction space has dimension
$(b-1)(g-d+r)$ as given by Lemma 5.3. In case $b=2$ and $\mu _0$ is injective,
the only possibilities are $\dim V =1$ corresponding to $\sigma _0 = \pi
_1^*\bold b+ \pi _2^*\bold d$, and $\dim V = 2$, $V = H^0(\scrO_{\Pee ^1}(1))$,
corresponding to $D_0 = D_1 + \pi _2^*\bold d'$, where $D_1$ is the
graph of a nonconstant map from $C$ to $\Pee ^1$. Thus the obstruction space
is necessarily zero in this case, provided that $\mu _0$ is injective.

Let us now describe all of the components of $H_{D_0, X}$ in the case
$b=2$. In this case
$0< -a < (g-1)/2$ and so $d > (g-1)/2$. Now suppose that $D_0 = D_1 + \pi
_2^*\bold d'$, where $D_1$ is the graph of a map from $C$ to $\Pee ^1$ of
degree
$d_1$ and $\bold d'$ has degree
$d' = d-d_1$. Thus $D_1$ corresponds to a linear subseries of $|\bold
d'|$, which we can write as $\Pee (V^*)$ for some vector space $V$ of dimension
two, together with a choice of isomorphism
$V \cong \Pee ^1$. In general we can think of $D_0$ as corresponding to a
sublinear system of $|\bold d|$ with base points. If
$C$ is generic in the Brill-Noether sense, then for there to exist a map from
$C$ to $\Pee ^1$ of degree $d_1$, we must have the Brill-Noether number $\rho =
2d_1 -g -2 \geq 0$, in which case the set
$G^1_d$ of all linear series of degree $d$ and dimension one has
dimension exactly $\rho$. Moreover for a generic curve $C$, if $\rho = 0$ then
$G^1_d$ consists of reduced points, whereas if $\rho >0$ then $G^1_d$ is smooth
and  irreducible of dimension  of dimension $\rho$ and the generic linear
series
in $G^1_d$ is complete. We see then that if $d = (g-1)/2, g/2, (g+1)/2$, then
$\frak M_0$ is the unique component of $H_{D_0, X}$. For $d= (g+2)/2$ the
components of $H_{D_0, X}$ are $\frak M_0$ together with a number of components
isomorphic to $\Pee ^3 = \Pee (H^0(\scrO_{\Pee ^1}(1) \otimes H^0(\bold d))$,
one for each $g^1_d$ on $C$. For $d> (g+2)/2$, there are two components,
$\frak M_0$ and a smooth component of dimension $\rho + 3 = 2d - g +1$, which
is essentially a $\Pee ^3$-bundle over  $G^1_d$.

Next we discuss the analytic structure of the moduli space in a neighborhood
of a singular point in case $b=2$, in other words how the components meet.
Recall
from Proposition 1.3 that the obstruction map has an intrinsically defined
quadratic term given by
$\alpha
\cup \partial
\alpha$.

\lemma{5.6} Let $D_0$ be given by the
section $\sigma _0 = \sum _i\alpha _i \otimes \beta _i\in H^0(\bold b) \otimes
H^0(\bold d)$, where the
$\alpha _i$ and $\beta _i$ are linearly independent and $V =
\operatorname{span}\{\alpha _i\}$ is a quotient of $H^0(\bold d)^*$ with dual
space $V^*\subseteq H^0(\bold d)$. Suppose that
$V \neq H^0(\bold d)$ and that the map
$\mu _0$ is injective for $\bold d$. Then cup product induces a surjective map
$$\gather
H^0(\scrO_X(D_0))/\Cee \sigma _0 \otimes \partial (H^0(\scrO_{D_0}(D_0))
\to H^1(\scrO_X(D_0))/\sigma _0H^1(\scrO_X) \\
\cong (H^0(\bold b)/V) \otimes H^1(\bold d).
\endgather$$
More precisely, for every $\tau _2 \notin V \subset H^0(\bold d)$,
the map
$$(H^0(\bold b)\otimes \Cee \tau _2) \otimes \partial (H^0(\scrO_{D_0}(D_0))
\to (H^0(\bold b)/V) \otimes H^1(\bold d)$$
is surjective.
\endstatement
\proof The image of $\partial \: H^0(\scrO_{D_0}(D_0)) \to H^1(\scrO_X) \cong
H^0(\scrO_{\Pee ^1})\otimes H^1(\scrO_C)$ is the set of all $1\otimes \xi$,
where $\sum _i\alpha \otimes \xi \beta _i =0$. Since the $\alpha _i$ are
linearly independent, this condition is equivalent to the condition that $\xi
\beta _i =0$ for all $i$. Cup product of such a class $1\otimes \xi$ with
$$\tau _1 \otimes \tau _2 \in H^0(\scrO_X(D_0))\cong H^0(\bold b) \otimes
H^0(\bold d)$$ is $\tau _1\otimes (\tau _2\cdot \xi) \in H^1(\scrO_X(D_0))$.
For the projection of this map to $$H^1(\scrO_X(D_0))/\sigma _0H^1(\scrO_X)
\cong  (H^0(\bold b)/V)
\otimes H^1(\bold d)$$ to be surjective, it suffices that, setting
$$(V^*)^\perp = \{\,\xi \in H^1(\scrO_C): \xi \cdot \beta _i = 0 \text{
for all $i$}\,\},$$ the induced cup product map $H^0(\bold d) \otimes
(V^*)^\perp \to H^1(\bold d)$ is surjective. Now by assumption, the adjoint
$\mu
_0^*$ of the
$\mu _0$ map is surjective, where
$$\mu _0^*\: H^1(\scrO_C) \to H^1(\bold d) \otimes H^0(\bold d)^*= \Hom
(H^0(\bold d), H^1(\bold d)).$$ Since $H^0(\bold d) \neq V$, there
exists a $\tau _2 \notin V$. Given $\eta \in H^1(\bold d)$, there
exists a linear map $F\: H^0(\bold d) \to H^1(\bold d)$ such that $F(V)
= 0$ and $F(\tau _2) = \eta$. The surjectivity of the map $\mu _0^*$ implies
that $F$ is given by taking cup product with a $\xi$ such that $\xi \cdot
\beta _i = 0$ for all $i$ and $\xi \cdot \tau _2 = \eta$. Thus the image of
$H^0(\scrO_X(D_0))/\Cee \sigma _0 \otimes \partial (H^0(\scrO_{D_0}(D_0))$
contains every element of the form $\tau _1 \otimes \eta$, where $\tau _1$ and
$\eta$ are arbitrary, and so is all of $(H^0(\bold b)/V)
\otimes H^1(\bold d)$.
\endproof

Using Lemma 5.6, we can describe the local structure of $H_{D_0, X}$ in case
$b=2$ near a reducible divisor $\pi _1^*\bold b +
\pi _2^*\bold d$, provided that $\bold d$ is generic in the sense that
$h^0(\bold d)$ is exactly $2$.

\corollary{5.7} With assumptions on the $\mu _0$ map as above, suppose that
$b=2$ and that $D_0 = \pi _1^*\bold b +
\pi _2^*\bold d$, where $h^0(\bold d) = 2$. Then an analytic neighborhood of
$H_{D_0, X}$ near $D_0$ is biholomorphic to a neighborhood of the origin in
$L_0
\cup L_1 \subset \Cee ^{d+2}$, where $L_0$ is a hyperplane and $L_1$ is a
linear space of dimension $2d-g +1$, not contained in $L_0$.
\endstatement
\proof The dimension of the Zariski tangent space $T$ to $H_{D_0, X}$ is
$2\cdot
2 -1 + d - 1 = d+2$. The dimension of $\Pee ^1 \times C_d$ is $d+1$, and so the
image $L_0$ of the tangent space to $\Pee ^1 \times C_d$ at $D_0$ has the
expected dimension of a hyperplane in $T$. In fact, we have seen in the
discussion prior to Lemma 5.2 that
$L_0$ is indeed a hyperplane, defined by the linear form
$\ell$, say. Thus if $\Phi\: T \to \Cee ^{g-d+1}$ is the Kuranishi obstruction
map, defined in a neighborhood of the origin, then there exists a
holomorphic function $f$ with differential $\ell$ such that $\Phi = f\Psi$, and
the quadratic term in $\Phi$ is equal to $\ell\cdot d\Psi _0$. The span
of $\alpha \cup \partial \alpha$ over all $\alpha$ is thus contained in the
image of $d\Psi _0$ and this span is the same, after polarizing, as the image
of
$\alpha \cup \partial\beta +\beta \cup \partial \alpha$ over all $\alpha,
\beta$. Using Lemma 5.6, there exists a choice of $\alpha _i, \beta _i$ with
$\partial \alpha _i = 0$ such that the obstruction space is generated by
$\alpha _i\cup \partial \beta_i$. Thus
$d\Psi _0$ has the same image as the map of Lemma 5.6 and so is surjective. It
follows that
$\Phi ^{-1}(0) = L_0 \cup \Psi ^{-1}(0)$, where $\Psi ^{-1}(0)$ is a smooth
submanifold of $T$ of codimension $g-d+1$. If it does not meet $L_0$
transversally, then $\Ker d\Psi _0 \subseteq L_0$. But $L_0
\cap H^0(\scrO_X(D_0))/\Cee \sigma _0$ is the tangent space to the Segre
embedding of $\Pee ^1 \times \Pee ^1$ in $\Pee H^0(\scrO_X(D_0)) = \Pee
^3$. Thus $H^0(\scrO_X(D_0))/\Cee \sigma _0$ is not contained in $L_0$. On
the other hand, $H^0(\scrO_X(D_0))/\Cee \sigma _0$ is the tangent space to
$\Pee H^0(\scrO_X(D_0))$ and so is unobstructed, so that it must be
contained in $\Ker d\Psi _0$. Thus $\Ker d\Psi _0$ is not contained in
$L_0$, and so $\Psi ^{-1}(0)$ meets $L_0$ transversally. This concludes the
proof.
\endproof

Finally, for a generic curve $C$, we shall use the description of the
components of the moduli space above to make some calculations in case $b=2$
and
$d$ is small. We do not need the description of the analytic structure of
the moduli space. In this case the moduli space always has the component
$\frak M_0= \Pee ^1\times \Sym ^dC$. For
$d= g-1/2, g/2, g+1/2$, $\frak M_0$ is the unique component, whereas in general
the moduli space is equal to $\frak M_0 \cup \frak M_1$. For $d= g+2/2$, $\frak
M_1$ is a union of $k$ copies of $\Pee ^3$, where $k$ is the number of
$g^1_d$'s on $C$, and for $d > g+2/2$ $\frak M_1$ is irreducible and
smooth, of the expected dimension $\rho +3 = 2d-g +1$.

To calulate the value of the Seiberg-Witten invariant, we shall first calculate
the contribution from $\frak M_0$ and then the contribution from $\frak M_1$.
Note that $\frak M_0$ does not have the expected dimension, which is $2d-g+1$.
Moreover the moduli space is in general singular. However it is easy to see
that we can choose incidence divisors $\mu _1,
\dots,
\mu _{2d-g+1}$ which meet properly in the smooth part of the moduli space.
Following the procedure of Section 3 (see the comments at the end of the
section), we first calculate the top Chern class of the obstruction bundle over
$\frak M_0$, which we have seen (Lemma 5.4) is the bundle $p_1^*T_{\Pee ^1}
\times p_2^*\Cal E$, at least after cutting down by
$\mu ^{2d-g+1}$. Here the
$p_i$ are the projections of $\Pee ^1\times \Sym ^dC$ to the $i^{\text{th}}$
factor and
$\Cal E$ is the bundle $R^1\psi _2{}_*\scrO_{C\times C_d}(\Cal I)$, of rank
$g-d$. The top Chern class of the tensor product of
$p_2^*\Cal E$ with the line bundle $p_1^*T_{\Pee ^1}$ is:
$c_N(p_1^*T_{\Pee ^1} \times p_2^*\Cal E) = \sum _{i=0}^Np_1^*c_1(T_{\Pee
^1})^ip_2^*c_{N-i}(\Cal E)$.
If $h$ is the hyperplane class on $\Pee
^1$, in other words the class of a point, then $p_1^*c_1(T_{\Pee ^1}) =
2p_1^*h$
and
$p_1^*c_1(T_{\Pee ^1})^i = 0$ for
$i>1$. Thus
$$c_N(p_1^*T_{\Pee ^1} \times p_2^*\Cal E) =p_2^*c_N(\Cal E) +
2p_1^*hp_2^*c_{N-1}(\Cal E).$$

By \cite{1}, p\. 322,
$c(\Cal E) = c(T_{C_d})^{-1} = (1+x)^{g-1-d}e^{\theta /1+x}$,
where $x$ is the class of the divisor $C_{d-1}\subset C_d$ and $\theta$ is the
pullback of the theta divisor on $\Pic^d C$ under the natural map, and
moreover $\theta ^kx^{d-k} = \frac{g!}{(g-k)!}$. To calculate the
Seiberg-Witten
invariant, we take $c_N(p_1^*T_{\Pee ^1} \times p_2^*\Cal E)$, where $N = g-d =
1-a$. This gives a class in $H^{2N}( \Pee ^1\times \Sym ^dC)$, and then we
further multiply by $\mu ^{d+1-N}$ and evaluate over the fundamental class. On
$\Pee ^1\times
\Sym ^dC$, it is clear that $\mu = p_1^*h+p_2^*x$, since $(t,p) \in \pi
_1^*\{s\} + \pi _2^*\bold d$ if and only if either $t=s$ or $p\in \bold d$, and
it is easy to see that the multiplicity is one. Thus we must calculate
$$\gather
(p_1^* h+p_2^*x)^{d+1-N}p_2^*c_N(\Cal E) +
2p_1^*hp_2^*c_{N-1}(\Cal E))\\
 = (2p_1^*h)p_2^*(c_{N-1}(\Cal E)x^{d+1-N} +
(d+1-N) p_1^*hp_2^*(c_N(\Cal E)x^{d-N}\\
= 2c_{N-1}(\Cal E)x^{d+1-N} + (d+1-N)c_N(\Cal E)x^{d-N}.
\endgather$$
Plugging in for $c(\Cal E)$, we have
$$\align
c(\Cal E) &= (1+x)^{g-1-d}e^{\theta /1+x} =(1+x)^{N-1}e^{\theta /1+x}\\
&= \sum_{i=0}^{N-1}\binom{N-1}{i}x^i\sum_{j=0}^\infty\frac{1}{j!}\theta
^j\sum_{k=0}^\infty\binom{-j}{k}x^k .
\endalign$$
Thus for example the term involving $c_{N-1}(\Cal E)$ becomes
$$\align
2&\sum _{i+j+k=
N-1}\binom{N-1}{i}\frac{1}{j!}\binom{-j}{k}x^{i+k+d+1-N}\theta ^j\\
=2&\sum _{i+j+k=
N-1}\binom{N-1}{i}\frac{1}{j!}\binom{-j}{k}x^{d-j}\theta ^j\\
=2&\sum _{i+j+k=
N-1}\binom{N-1}{i}\binom{-j}{k}\frac{1}{j!}\frac{g!}{(g-j)!}\\
=2&\sum _{j=0}^{N-1}\left(\sum
_{k=0}^{N-1-j}\binom{N-1}{N-1-j-k}\binom{-j}{k}\right)\binom{N+d}{j},
\endalign$$
where we have used $g=N+d$.

Applying Lemma 4.5  with $e=0$ to the inner sum above, where we let $a =
N-1-j$ for a fixed $j$, we see that the expression reduces to $\dsize 2\sum
_{j=0}^{N-1}\binom{N+d}{j}$.

A very similar manipulation with the term $(d+1-N)c_N(\Cal E)x^{d-N}$ gives
$$
(d+1-N)c_N(\Cal E)x^{d-N} =
 (d+1-N)\binom{N+d}{N}.$$
 The final contribution for the component $\frak M_0$ is
therefore
$$2\sum _{j=0}^{N-1}\binom{N+d}{j}+ (d+1-N)\binom{N+d}{N}.$$
Note that $\frak M_0$ is the unique component for $N= d+1, d, d-1$
corresponding to the cases $g = 2d+1, 2d, 2d-1$. Plugging in, we find that the
value of the invariant in case $N= d+1$ is
$$2\sum _{j=0}^{d}\binom{2d+1}{j}= \sum _{j=0}^{2d+1}\binom{2d+1}{j} =
(1+1)^{2d+1} = 2^g.$$
Similar calculations handle the cases $N= d, d-1$, and again give the value
$2^g$.

For $d= g+2/2,N=d-2, g = 2d-2$, the set $\frak M_1$ consists of $k$ copies of
$\Pee ^3$, where
$k$ is the number of $g^1_d$'s on $C$. This number has been computed by
Castelnuovo \cite{1} p\. 211: it is
$$g!\frac1{(g-d+1)!}\frac1{(g-d+2)!} = \frac{(2d-2)!}{(d-1)!\,d!}.$$
In this case, the restriction of $\Cal D$ to each piece $X\times \Pee ^3$ is
the incidence divisor, so that $\mu$ restricts to the hyperplane class in each
$\Pee ^3$. Thus the final answer is
$$2\sum _{j=0}^{d-3}\binom{2d-2}{j} +3\binom{2d-2}{d-2}+
\frac{(2d-2)!}{(d-1)!\,d!}$$
which after a brief manipulation becomes
$\dsize\sum _{j=0}^{2d-2}\binom{2d-2}{j} = (1+1)^{2d-2} =
2^g$.

Somewhat more involved methods handle the case $N=d-3$, and presumably might be
pushed, using excess intersections, to give the general case. However, we shall
give a simpler method for the calculation in the next section.

\section{6. Deformation to more general ruled surfaces.}

In this section, we shall study Seiberg-Witten moduli spaces, or equivalently
the Hilbert scheme, for more general ruled surfaces $\Pee (V)$, where $V$ is a
general (and in particular stable) rank two bundle over $C$. We shall
deal with the case where $\det V$ has even degree, and thus assume that $c_1(V)
= 0$.  Also, we shall only
discuss the case of sections of $V$. However, it will be clear that our
methods generalize to handle the case of odd degree as well as more general
cases of multisections, and thus suffice for the homological calculations of
the
invariants in general. We will outline this approach at the end.
Throughout this section, we fix a smooth curve $C$ of genus $g\geq 2$. it will
not be necessary to assume that $C$ is generic in the Brill-Noether sense.

Recall that there is a one-to-one correspondence between irreducible sections
$D_0$ of $\Pee (V)$ and line bundles $\lambda$ such that $V\otimes \lambda$ has
a nowhere vanishing section, as follows: given a section $D_0$ of $\Pee (V)$,
apply $R^i\pi _*$ to the exact sequence
$$0 \to \scrO_X \to \scrO_X(D_0) \to \scrO_{D_0}(D_0) \to 0$$
to obtain the exact sequence
$$0 \to \scrO_C \to R^0\pi _*\scrO_X(D_0) \to\scrO_{D_0}(D_0) \to 0.$$
Here $R^1\pi _*\scrO_X= 0$ since the fibers are $\Pee ^1$ and we can write
$R^0\pi _*\scrO_X(D_0) = V\spcheck \otimes \lambda = V\otimes \lambda$. Note
that the normal bundle $\scrO_{D_0}(D_0)$ is naturally $\lambda ^2$. The
inverse map sends the section of $V\spcheck \otimes
\lambda$ to the homogeneous degree one subvariety of $\Pee(V)$ that it defines.

If
$D_0$ is not irreducible, then $D_0 = E_0 + \pi ^*\bold e$. In this case the
map
$\scrO_X \to \scrO_X(D_0)$ factors through $\scrO_X(E_0)$ and the induced map
$\scrO_C \to R^0\pi _*\scrO_X(D_0)$ factors as
$$\scrO_C \to R^0\pi _*\scrO_X(E_0)\to R^0\pi _*\scrO_X(E_0)\otimes
\scrO_C(\bold e) = R^0\pi _*\scrO_X(D_0).$$
Thus $D_0$ still corresponds to a section of $V\spcheck\otimes \lambda$ for an
appropriate $\lambda$, but the section vanishes exactly along $\bold e$. Again
such a section defines a subvariety of $\Pee (V)$, which is exactly $D_0$ since
the section vanishes along $\bold e$. In this way, we can identify $|D_0|$ with
$\Pee H^0(V\spcheck \otimes \lambda)$, including the reducible fibers.

\proposition{6.1} Let $e$ be an positive integer. For every line bundle
$\lambda$ on $C$ of degree $e$,
\roster
\item"{(i)}" There exist stable bundles $V$ on $C$ together with an exact
sequence
$$0 \to \lambda ^{-1} \to V \to \lambda \to 0.$$
\item"{(ii)}" For $e< (g-1)/2$, and $V$ a generic stable bundle satisfying
\rom{(i)}, if $\mu$ is a line bundle of degree $d\leq e$ and $H^0(V\otimes \mu)
\neq 0$, then $\mu = \lambda$ and $H^0(V\otimes \lambda)$ has dimension one.
\item"{(iii)}" For $e = (g-1)/2$ and $V$ general, there are exactly $2^g$
distinct $\lambda$ with  $H^0(V\otimes \lambda) \neq 0$, and for each such
$\lambda$,
$\dim H^0(V\otimes \lambda) =1$. Moreover, if $\deg \mu < (g-1)/2$, then
$H^0(V\otimes \mu) = 0$.
\endroster
\endstatement
\proof For a line bundle $\lambda$ of degree $e>0$, $\Ext ^1(\lambda, \lambda
^{-1}) = H^1(\lambda ^{-2})$ which has dimension $2e+g-1$, by Riemann-Roch. Let
$V$ be a rank two bundle corresponding to an extension class $\xi \in
H^1(\lambda ^{-2})$. Suppose that there exists a nonzero map $\mu \to V$, where
$\deg \mu = d\leq e$, and that $\mu \neq \lambda$ in case $d=e$. Since $H^0(\mu
\otimes \lambda ^{-1})=0$ since $\mu \otimes \lambda ^{-1}$ either has negative
degree or has degree zero and is not trivial, there must exist a nonzero
section
of $\mu \otimes \lambda$ which lifts to a section of $\mu \otimes V$. Note that
$(s)$, the divisor of zeroes of $s$, has degree $d+e$, and we can identify the
set of all pairs $(\mu, s)$ such that $\mu$ is a line bundle of degree $d$ and
$s$ is a nonzero section of $\mu \otimes \lambda$, mod scalars, with $C_{d+e}$.

The section $s$ lifts to a section of $\mu \otimes V$ if and only if the
coboundary map $\partial (s) = 0$, where $\partial s \in H^1(\mu \otimes
\lambda ^{-1})$. Now $\partial s = \xi \cdot s$, the cup product of $s\in
H^0(\mu \otimes \lambda)$ with $\xi \in H^1(\lambda ^{-2})$. Consider the exact
sequence
$$0 \to \lambda ^{-2} @>{\times s}>> \lambda ^{-1}\otimes \mu \to \scrO_{\bold
f} \to 0,$$
where $\bold f = (s) \in C_{d+e}$. By assumption, $\partial s = 0$ if and only
if $\xi \cdot s = 0$ if and only if the image of $\xi$ in $H^1(\lambda
^{-1}\otimes \mu)$ is zero, if and only if $\xi$ is in the image of
$H^0(\scrO_{\bold f})$. By assumption either $\deg (\lambda ^{-1}\otimes \mu) <
0$ or $\deg (\lambda ^{-1}\otimes \mu)  = 0$ but $\lambda ^{-1}\otimes \mu$ is
not trivial. In either case $H^0(\lambda ^{-1}\otimes \mu) = 0$, and so the
image of $H^0(\scrO_{\bold f})$ has dimension $d+e$. Thus the set of possible
extension classes $\xi$ for which a given
$s$ lifts has dimension $d+e$, and so corresponds to a $\Pee ^{d+e-1} \subseteq
\Pee H^1(\lambda ^{-2})= \Pee ^{2e+g-2}$. The set of all $\xi$ for which some
$s$ lifts is then the union over all $s$ of a linear subspace of $\Pee
^{2e+g-2}$ of dimension $d+e-1$. Since the set of all $s$ is just $C_{d+e}$,
the dimension of the set of all possible $\xi$ is at most $2d+2e-1$. This
number is  less than $2e+g-2$ exactly when $d < (g-1)/2$. Choosing a bundle $V$
coresponding to an extension class $\xi$ in the complement of this set  gives a
bundle $V$, written as an extension of $\lambda$ by $\lambda ^{-1}$, such that,
if $H^0(V\otimes \mu)
\neq 0$ and $\deg \mu \leq \lambda$, then $\mu = \lambda$. In particular, $V$
is stable, proving (i) and (ii), except for the statement that $\dim H^0
(V\otimes \lambda ) = 1$.

To see the statement about $\dim H^0(V\otimes \lambda)$,  if $s$ is a nonzero
section of $\lambda ^2$ which lifts to a section of $V\otimes \lambda$, then
arguments similar to those above show that the orthogonal complement of $s\cdot
H^0(K_C)$ in $H^1(\lambda ^2)$ has dimension $2e-1$ and so gives a linear
space of dimension $2e-2$ inside $\Pee ^{2e+g-2}$. Moreover, the possible $s$
correspond to the case $\mu = \lambda$, and so form a proper subvariety of
$C_{d+e} = C_{2e}$. In all, the $\xi$  for which some $s$ lifts, such that the
corresponding line bundle $\mu = \lambda$, form a subvariety of $\Pee
^{2e+g-2}$ of dimension at most $4e-3$. Now $4e-3< 2e+g-2$ provided that
$e<(g+1)/2$, and thus certainly if $e\leq (g-1)/2$. This establishes the last
statement of (ii).

The remaining assertion (iii) is a classical formula due to Corrado Segre
\cite{11}. It follows from our
formulas in the previous section for the zero-dimensional invariant, and will
be reproved in more generality shortly.
\endproof

We now fix a bundle $V$ which will later be assumed generic in an appropriate
sense. For each $d$, we consider the following varieties of Brill-Noether type:
$$\align
W_{1,d}(V) &= \{\, \lambda \in \Pic ^d C: h^0(V\otimes \lambda ) \geq 1\,\};\\
G_{1,d}(V) &= \{\,(s,\lambda): \lambda \in \Pic ^d C, s\in \Pee(H^0(V\otimes
\lambda ))\,\}.
\endalign$$
Thus there is a natural map $G_{1,d}(V) \to W_{1,d}(V)$. It is also clear
that, with $X= \Pee (V)$, $G_{1,d}(V)$ is exactly the Hilbert scheme of $X$
corresponding to sections (possibly reducible) of the appropriate degree
and that the map $G_{1,d}(V) \to W_{1,d}(V)$ can be identified with the map
from the Hilbert scheme to $\Pic X$.   We wish to give another construction of
the Hilbert scheme in this context; in other words, we will put another scheme
structure on $G_{1,d}(V)$ and then claim that it is in fact the usual one. To
do so, we make a construction similar to the usual construction of
Brill-Noether theory: fix a divisor $D$ on $C$ of degree $m\gg 0$ such that
$h^1(V\otimes \lambda \otimes \scrO_C(D)) = 0$ for all line bundles $\lambda$
of
degree $d$. We can assume that $D$ is an effective divisor consisting of
reduced points of $C$ if we wish. Consider the restriction sequence
$$0 \to V\otimes \lambda \to V\otimes \lambda \otimes \scrO_C(D) \to V\otimes
\lambda \otimes \scrO_D(D)\to 0.$$
Taking global sections, there is a map
$$\varphi\: H^0(V\otimes \lambda \otimes \scrO_C(D)) \to H^0( V\otimes
\lambda \otimes \scrO_D(D)).$$
The first vector space has dimension $2m+2d-2g+2$, by Riemann-Roch, and the
second has dimension $2m$, and $\lambda \in W_{1,d}(V)$ if and only if
$\varphi$ has a kernel. In this case, the fiber over $\lambda$ in
$G_{1,d}(V)$ is just $\Pee (\Ker \varphi) = \Pee(H^0(V\otimes \lambda))$.
Globally, let
$\Cal P$ be a Poincar\'e line bundle for
$C\times \Pic ^dC$, and let $\pi _i$ be the projection of $C\times \Pic ^dC$ to
the $i^{\text{th}}$ factor. Set
$$\align
\Cal E' &= \pi _2{}_*\left(\Cal P\otimes \pi _1^*(V\otimes
\scrO_C(D))\right);\\
\Cal E'' &= \pi _2{}_*\left(\Cal P\otimes \pi _1^*(V\otimes \scrO_D(D))\right),
\endalign$$
so that there is a natural evaluation map $\Phi \: \Cal E' \to \Cal E''$. Then
$W_{1,d}(V)$ is the scheme where $\Phi$ fails to be injective. We can define
$G_{1,d}(V)$ similarly: let $p\: \Pee\Cal E' \to \Pic^dC$ be the projection. We
have the inclusion of $\scrO_{\Pee\Cal E'}(-1)$ inside $p^*\Cal E'$. Consider
the composition
$$\scrO_{\Pee\Cal E'}(-1) \to p^*\Cal E' @>{p^*\Phi}>>p^*\Cal E''.$$
If we denote this composition by $\tilde \Phi$, then $\tilde \Phi =0$ at a
point $(s, \lambda)$, where $\lambda \in \Pic ^dC$ and $s \in \Pee H^0(V\otimes
\lambda \otimes \scrO_C(D))$, if and only if $s$ is in the image of
$\Pee H^0(V\otimes \lambda)$. Thus the vanishing of $\tilde \Phi$ defines
$G_{1,d}(V)$ as a set inside $\Pee \Cal E'$. Note that $\Pee \Cal E'$ is
itself a Hilbert scheme: it is the same as $G_{1, d+m}(V)$, corresponding to
the set of all sections of $X$ of degree $2d+2m$. Moreover $G_{1,d}(V)$ is the
subset of $\Pee \Cal E'$ consisting exactly of those sections containing $\pi
^*D$. We leave it to the reader to work through the details that the subscheme
defined by
$\tilde
\Phi$ represents the functor corresponding to $G_{1,d}(V)$ (see \cite{1} pp\.
182--184 for the Brill-Noether analogue) and that indeed this identifies
$G_{1,d}(V)$ with the Hilbert scheme as schemes.

Next suppose that $(s, \lambda)$ is a point of $G_{1,d}(V)$ such that
the section $s$ does not vanish. Standard arguments (cf\.
\cite{1}, pp\. 185--186) identify the Zariski tangent space to
$G_{1,d}(V)$ at the point $(s, \lambda)$ with $H^0(\lambda ^2)= \Hom (\lambda
^{-1}, V/\lambda ^{-1})$, via the exact sequence
$$0\to H^0(\scrO_C) \to H^0(V\otimes \lambda) \to H^0(\lambda ^2) \to
H^1(\scrO_C) \to H^1(V\otimes \lambda).$$
Moreover the differential of the map from $G_{1,d}(V)$ to $W_{1,d}(V)$ is the
obvious map $ H^0(\lambda ^2) \to H^1(\scrO_C)$. Finally,  a standard
cocycle calculation identifies the obstruction space as $H^1(\lambda ^2)$. Note
that, if $\lambda$ corresponds to the section $D_0$ of $X$, then  $H^1(\lambda
^2) = H^1(\scrO_{D_0}(D_0))$ is the same obstruction we would have found via
the Hilbert scheme, as well it must be since $G_{1,d}(V)$ represents the same
functor as the Hilbert scheme.

For a general section
$s$, suppose that the map $\scrO_C \to V\otimes \lambda$ vanishes along the
divisor
$\bold e$, so that there is a factorization
$$\scrO_C \to V\otimes \lambda \otimes \scrO_C(-\bold e) \to V\otimes
\lambda.$$
Let $\lambda _0 = \lambda \otimes \scrO_C(-\bold e)$. Then there is a
commutative diagram
$$\CD
0 @>>> \scrO_C @>>> V\otimes \lambda \otimes \scrO_C(-\bold e) @>>> \lambda
_0^2 @>>> 0\\
@. @| @VVV @VVV @.\\
0 @>>> \scrO_C @>>> V\otimes \lambda @>>> \lambda _0^2\oplus T @>>> 0.
\endCD$$
Here $T$ is a skyscraper sheaf isomorphic to $V\otimes \scrO_{\bold e}$, and
so has length $2\deg \bold e = 2e$, say. We can again identify the Zariski
tangent space to $G_{1,d}(V)$ at  $(s, \lambda)$ with
$$\Hom (\lambda ^{-1}, V/\lambda ^{-1}) = H^0 (V\otimes \lambda/s\cdot
\scrO_C)=H^0(\lambda_0 ^2\oplus T) = H^0(\lambda_0 ^2)
\oplus H^0(T).$$
Moreover the obstruction space is
$$\Ext ^1(\lambda ^{-1}, V/\lambda ^{-1}) = H^1(\lambda_0 ^2\oplus T) =
H^1(\lambda_0 ^2).$$ This again corresponds to the deformation
theory and obstruction theory for the Hilbert scheme: let $D_0 = E_0 + \pi
^*\bold e$, where $E_0$ is an irreducible section of $X$. Apply $R^i\pi _*$ to
the exact sequence
$$0 \to \scrO_X \to \scrO_X(D_0) \to \scrO_{D_0}(D_0) \to 0,$$
using $R^1\pi _*\scrO_X = R^1\pi _*\scrO_X(D_0) = 0$. We obtain
$$0 \to \scrO_C \to  V\otimes \lambda \to \pi _*\scrO_{D_0}(D_0) \to 0,$$
and $R^1\pi _*\scrO_{D_0}(D_0)  = 0$. Thus $H^i(\scrO_{D_0}(D_0) ) = H^i(\pi
_*\scrO_{D_0}(D_0) )= H^i(V\otimes \lambda/s\cdot \scrO_C)$.

Summarizing, we have shown the following:

\proposition{6.2} Let $(s, \lambda)$ be a point of $G_{1,d}(V)$. Suppose that
the section $s$ vanishes exactly along the effective divisor $\bold e$, and set
$\lambda _0 = \lambda \otimes \scrO_C(-\bold e)$. Then the Zariski tangent
space to $G_{1,d}(V)$ is $H^0(V\otimes \lambda/s\cdot \scrO_C)$, which has
dimension $h^0(\lambda _0^2) + 2e$, and the obstruction space is $H^1(\lambda
_0^2)$. \qed
\endstatement

Assuming for simplicity that $s$ does not vanish at any point, and so $\lambda
= \lambda _0$ in the above notation, the group $H^1(\lambda ^2)$ arises in yet
another way as follows: the universal extension over $\Pee H^1(\lambda ^2)$ of
$\lambda$ by $\lambda ^{-1}$ gives rise to a Kodaira-Spencer map from the
tangent space of $\Pee H^1(\lambda ^{-2})$ at a nonzero point $\xi \in
H^1(\lambda ^{-2})$, namely  $H^1(\lambda ^{-2})/\Cee \cdot \xi$, to
$H^1(\operatorname{ad}V)$. A diagram chase identifies the cokernel of this map
with $H^1(\lambda ^2)$ under the natural map $H^1(\operatorname{ad}V) \to
H^1(\lambda ^2)$. Thus $H^1(\lambda ^2)=0$ if and only if the map from the set
of extensions to moduli is a submersion at $\xi$. If $d\geq (g-1)/2$, then for
a generic choice of $\lambda$ we will indeed have $H^1(\lambda ^2)=0$, and so
the map from extensions to moduli will be a submersion where defined.

We now show that, for a generic choice of $V$, the Hilbert scheme is always
smooth:

\proposition{6.3} For a generic stable bundle $V$ with $c_1(V) = 0$ and for
all $d$, the Hilbert scheme of sections of $X= \Pee(V)$ of square $2d$ is
everywhere smooth of the expected dimension $2d-g+1$.
\endstatement
\proof Fix a value for $d$ corresponding to the degree of a line subbundle
$\lambda$ of $V$. We have seen that, for generic
$V$, the space of sections is empty if $d< (g-1)/2$, and if $\lambda$ has
degree $>g-1$, then $h^1(\lambda ^2) = 0$ by Serre duality. For $(g-1)/2 \leq
d \leq g-1$, let $\Cal P$ be the Poincar\'e line bundle over $C\times \Pic^dC$
and let $\Cal V = R^1\pi _2{}_*\Cal P^{\otimes -2}$. Note that, as $\deg
\lambda = d >0$, then $h^0(\lambda ^{\otimes -2}) = 0$ and $h^1(\lambda
^{\otimes -2}) =  2d+g-1$. Hence $\Cal V$ is a vector bundle of rank $2d+g-1$
over $\Pic^dC$ and $\Pee\Cal V$  is a
$\Pee ^{2d+g-2}$-bundle over $\Pic^dC$. The space $\Pee\Cal V$ is a moduli
space for vector bundles $V$ given as extensions. There is an open
subset
$\Cal U$ of $\Pee\Cal V$ corresponding to stable bundles, which is nonempty
by (i) of (6.1). The remarks prior to the statement of (6.3) imply that the
morphism $\Cal U \to \frak M(C)$ is dominant, where
$\frak M(C)$ is the moduli space of stable rank two bundles $V$ over $C$ with
$c_1(V) = 0$. Let
$$\Cal B = \{\, \lambda \in \Pic^dC: h^1(\lambda ^2) \neq 0\,\}.$$
Since $h^1(\lambda ^2) \neq 0$ if and only if $h^0(K_C\otimes \lambda ^{-2})
\neq 0$, $\Cal B$ is the inverse image in $\Pic^dC$ of the set of effective
divisors in $\Pic^{2g-2-2d}C$ under the obvious (\'etale) map $\lambda \mapsto
K_C\otimes \lambda ^{-2}$. Thus $\Cal B$ has the same dimension as the set of
effective divisors in $\Pic^{2g-2-2d}C$, namely $2g-2-2d$. Hence if $p\: \Pee
\Cal V \to \Pic ^dC$ is the projection, then
$$\dim p^{-1}(\Cal B) = 2g-2 -2d + 2d+g-2 =3g-4.$$
It follows that the image of $p^{-1}(\Cal B)\cap \Cal U$ in $\frak M(C)$ cannot
be all of $\frak M(C)$. Thus we can choose a stable bundle $V$ such that, if
$\deg
\lambda <(g-1)/2$, then there is no nonzero section of $V\otimes
\lambda$, and if $\deg \lambda \geq (g-1)/2$ and
there is a nowhere vanishing section $s$ of $V\otimes \lambda$, then $
h^1(\lambda ^2) = 0$. It follows that $(s, \lambda)$ is a smooth point of
$G_{1,d}(V)$ (or of the Hilbert scheme), and the discussion prior to
Proposition
6.2 shows how to extend this to all nonzero sections. Thus the Hilbert scheme
of sections is everywhere smooth.
\endproof

Next we turn to the enumerative geometry of the Hilbert scheme. Since
$\scrO_{\Pee\Cal E'}(-1)$ is a line bundle, $\tilde \Phi$ is equivalent to a
section of
$\scrO_{\Pee\Cal E'}(1)\otimes \Cal E''$ and the class of its zero set, namely
$G_{1,d}(V)$, is given by $c_{2m}(\scrO_{\Pee\Cal E'}(1)\otimes \Cal E'')$. If
we set $\zeta = c_1(\scrO_{\Pee\Cal E'}(1))$ and use the fact that $\Cal E''$
is
a topologically trivial bundle of rank $2m$ \cite{1} p\. 309, then
$c_{2m}(\scrO_{\Pee\Cal E'}(1)\otimes \Cal E'')$ is the term of degree $2m$ in
$(1+\zeta)^{2m}$, namely $\zeta ^{2m}$. In particular, the class of
$W_{1,d}(V)$ is given by $p_*\zeta ^{2m}$, which is the appropriate Segre class
of $\Cal E'$. To find it, take $c(\Cal E')^{-1}$. Topologically $\Cal E'$ is
two copies of $\pi _2{}_*\left(\Cal P \otimes \pi _1^*\scrO_C(D)\right)$, and
using \cite{1} p\. 336, this last bundle has total Chern class $e^{-\theta}$.
Thus
$c(\Cal E') = e^{-2\theta}$ and $c(\Cal E')^{-1} = e^{2\theta}$. Taking the
term of degree $2g-2d-1$, we find \cite{6}:
$$[W_{1,d}(V)] = \frac{(2\theta)^{2g-2d-1}}{(2g-2d-1)!}.$$
In particular, when $2d+1 = g$ we obtain $(2\theta)^g/g! = 2^g$, giving the
formula of Segre for the case of the zero-dimensional invariant. (In fact,
Segre stated the formula for all values of $d$.)

To handle the general case, we use:

\proposition{6.4} The $\mu$ divisor on $G_{1,d}(V)$ is algebraically equivalent
to the restriction of $\zeta = c_1(\scrO_{\Pee\Cal E'}(1))$.
\endstatement

\corollary{6.5} The value of the Seiberg-Witten invariant is $2^g$.
\endstatement
\demo{Proof of the corollary} We need to compute $\zeta ^{2m+ 2d-2g+1}$. This
is the top Segre class of $\Cal E'$, and by the calculations above it is equal
to the degree $g$ term in $e^{2\theta}$, namely $2^g$.
\endproof

\demo{Proof of \rom{(6.4)}} Keeping our previous notation, note that $\Pee \Cal
E'$ is itself a Hilbert scheme, namely the scheme of all sections of $X$ of
degree $2d+2m$, and the Hilbert scheme is the subscheme of all sections which
are of the form $D_0 + \pi ^*D$. Thus, choosing a point $p\in X$ not lying in
$\pi ^*D$, the incidence divisor $\mu (p)$ for the Hilbert scheme is the
restriction of the corresponding incidence divisor on $\Pee \Cal E'$. Recall
that $\Cal E ' = \pi _2{}_*(\Cal P \otimes \pi _1^*(V\otimes \scrO_C(D)))$. Let
$t = \pi (p)$, and suppose that we have chosen the Poincar\'e line bundle
$\Cal P$ so that $\Cal P|\{t\}\times \Pic ^dC$ is trivial. Fix the line $\ell
\subset V_t$ corresponding to $p\in X$. Identifying the space $(V_t/\ell)
\otimes _\Cee  \scrO_C(D) _t$ with $\Cee$, there is then a surjection
$$\Cal P \otimes \pi _1^*(V\otimes \scrO_C(D)) \to \Cal P \otimes \pi
_1^*(V\otimes \scrO_C(D))|\{t\}\times \Pic ^dC \to \Cal P |\{t\}\times \Pic ^dC
\cong \scrO_{\Pic ^dC}.$$
Applying $\pi _2{}_*$, we get a map $F\: \Cal E' \to \scrO_{\Pic ^dC}$, which
is nonzero if $D$ is sufficiently ample. Clearly $F(s, \lambda) =0$ exactly
when $s(t) \in \ell \subset V_t$. In other words, the zero set of $F$ is the
incidence divisor $\Cal D$ corresponding to the point $p$. Now
$$\align
F \in &\Hom (\Cal E', \scrO_{\Pic ^dC}) = H^0(\Pic ^dC; (\Cal
E')\spcheck) \\
&= H^0(\Pic ^dC; p_*\scrO_{\Pee \Cal E'}(1)) = H^0(\Pee \Cal E'; \scrO_{\Pee
\Cal E'}(1)).
\endalign$$
Running through the identifications above, we see that the zero set of $F$ is
exactly the zero set of the induced section of $\scrO_{\Pee
\Cal E'}(1)$, and a straightforward argument also checks the multiplicity. Thus
$\Cal D$ is the zero set of a section of
$\scrO_{\Pee
\Cal E'}(1)$, and so $[\Cal D] = \zeta$.
\endproof

We can apply the methods above to handle enumerative
questions not directly related to Seiberg-Witten theory. For example, there is
the following calculation via the Grothendieck-Riemann-Roch theorem:

\theorem{6.6} Suppose that $g$ is even. Then the set of stable rank two bundles
$V$ with
$\det V=0$ such that there exists a line bundle $\lambda$ of degree $(g-2)/2$
with $h^0(V\otimes \lambda ) \neq 0$ is an irreducible divisor in the moduli
space $\frak M(C)$. Its class, at least as a divisor on the moduli functor,
is $2^g\Delta$, where $\Delta$ is the first Chern class of the determinant line
bundle, which again exists on the moduli functor.
\endstatement

Finally let us give the general formula for the Seiberg-Witten invariant:

\theorem{6.7} Let $X = \Pee ^1 \times C$ be a product ruled surface and let $L$
be a line bundle on $X$ of type $(2a,2b)$, with $b\geq 1$ and $(1-g)/b \leq a
<0$. Then for a K\"ahler metric $\omega$ such that $\omega \cdot L < 0$, the
value of the Seiberg-Witten invariant on $L$ is $b^g$.
\endstatement

Note that this formula is a special case of a general transition formula for
Seiberg-Witten invariants, which has been established by Li and Liu \cite{9} as
well as the authors (unpublished). It would also follow, by copying the
arguments
above for the case of sections, but working with $\Sym ^mV$ (for $m=
b-1$) instead of
$V$, if we knew that, for a general ruled surface, the Hilbert scheme was
always
smooth of the expected dimension. We state this as a conjecture:

\medskip
\noindent {\bf Conjecture.} Let $X$ be a general ruled surface. Then every
component of the Hilbert scheme is smooth of the expected dimension. Here,
since $H^2(\scrO_X) = 0$, the expected dimension of the Hilbert scheme at a
curve $D$ is $\dsize \frac12(D^2 - D\cdot K_X)$.
\medskip

Without assuming this conjecture, one can deduce the result from the methods of
Section 3 in case the Hilbert scheme is smooth but does not have the expected
dimension, and with more work in general, again by reducing it to a homological
calculation in a space along the lines of $\Pee \Cal E'$.


\Refs

\ref \no  1\by E. Arbarello, M. Cornalba, P. A. Griffiths and J. Harris \book
Geometry of Algebraic Curves volume I \publ Springer Verlag \publaddr New
York Berlin Heidelberg Tokyo \yr 1985 \endref

\ref \no 2 \by R. Brussee \paper Some $C^\infty$ properties of K\"ahler
surfaces \paperinfo Algebraic geometry e-prints 9503004 \endref

\ref \no 3\by R. Fintushel, P. Kronheimer, T. Mrowka, R. Stern, and C. Taubes
\toappear
\endref

\ref \no 4\by R. Friedman and J W. Morgan\book Smooth Four-Manifolds and
Complex Surfaces, {\rm Ergebnisse der Mathematik und ihrer Grenz\-gebiete 3.
Folge} {\bf 27} \publ Springer \publaddr Berlin Heidelberg  New York \yr
1994\endref

\ref \no 5\bysame \paper Algebraic surfaces and
Seiberg-Witten invariants \toappear \endref

\ref \no 6 \by F. Ghione \paper Quelques r\'esultats de Corrado Segre sur les
surfaces r\'egl\'ees \jour Math. Annalen \vol 255 \yr 1981 \pages 77--95
\endref

\ref \no 7\by K. Kodaira and D.C. Spencer \paper A theorem of completeness of
characteristic systems of complete continuous systems \jour Amer. J. Math \vol
81 \yr 1959 \pages 477--500 \endref

\ref \no 8 \by H. Lange \paper H\"ohere Sekantenvariet\"aten und Vektorb\"undel
auf Kurven \jour Manuscripta Math. \vol 52 \yr 1985 \pages 63--80 \endref

\ref \no 9 \by T. J. Li and A. Liu \paper General wall crossing formula
\paperinfo preprint \endref

\ref \no 10\by D. Mumford \book Lectures on Curves on an Algebraic Surface
\bookinfo Annals of Mathematics Studies \vol 59\publ Princeton University Press
\publaddr Princeton, NJ
\yr 1966
\endref

\ref \no 11 \by C. Segre \paper Recherches g\'en\'erales sur les courbes et les
surfaces r\'egl\'ees alg\'ebriques II \jour Math. Annalen \vol 34 \yr 1889
\pages 1--25 \endref

\ref \no 12\by E. Witten \paper Monopoles and four-manifolds \jour Math.
Research Letters \vol 1 \yr 1994 \pages 769--796
\endref


\endRefs
\enddocument